\documentclass[11pt]{article}
\usepackage{amssymb}
\usepackage{bbm}
\newtheorem{pr}{Proposition}
\newtheorem{de}{Definition}
\newtheorem{lem}{Lemma}
\newtheorem{them}{Theorem}
\newcommand{\nc}[2]{\newcommand{#1}{#2}}

\nc{\bsa}{\begin{satz}}
\nc{\bpr}{\begin{pr}}
\nc{\bth}{\begin{them}}
\nc{\ble}{\begin{lem}}
\nc{\bco}{\begin{corollary}}
\nc{\bre}{\begin{remark}}
\nc{\bex}{\begin{example}}
\nc{\bde}{\begin{de}}
\nc{\ede}{\end{de}}
\nc{\esa}{\end{satz}}
\nc{\epr}{\end{pr}}
\nc{\ethe}{\end{them}}
\nc{\ele}{\end{lem}}
\nc{\eco}{\end{corollary}}
\nc{\ere}{\hfill\mbox{$\Diamond$}\end{remark}}
\nc{\eex}{\end{example}}
\nc{\epf}{\hfill\mbox{$\square$}}
\nc{\beq}{\begin{equation}}
\nc{\eeq}{\end{equation}}
\nc{\ot}{\otimes}
\nc{\lra}{\longrightarrow}
\nc{\ci}{\circ}
\def\Z{\mathbbm{Z}}
\begin{document}
\title{Covering and gluing of algebras and differential algebras}
\author{Dirk Calow\thanks{supported by Deutsche Forschungsgemeinschaft, e-mail
Dirk.Calow@itp.uni-leipzig.de}
 \, and Rainer Matthes\thanks{supported by S\"achsisches Staatsministerium
 f\"ur Wissenschaft und Kunst,\newline\hspace*{.5cm} e-mail 
 Rainer.Matthes@itp.uni-leipzig.de}\\[.5cm]
\normalsize Institut f\"ur Theoretische Physik\\
\normalsize Fakult\"at f\"ur Physik und Geowissenschaften\\
\normalsize Universit\"at Leipzig\\
\normalsize Augustusplatz 10/11\\
\normalsize D-04109 Leipzig\\
\normalsize Germany}
\date{}
\maketitle
\begin{abstract}
Extending work of Budzy\'nski and Kondracki, we investigate coverings and gluings of
algebras
and differential algebras. We describe in detail the gluing of two quantum
discs along their classical subspace, giving a $C^*$-algebra isomorphic to a
certain Podle\'s sphere, as well as the gluing of
$U_{q^{1/2}}(sl_2)$-covariant
differential calculi on the discs.
\end{abstract}
{\small 
}
\section{Introduction}
\indent

In classical geometry or topology, the objects of study can always be
considered as being glued together from local pieces which are in a certain
sense (e.g. topologically) trivial. The extension of this idea to
noncommutative geometry seems to be not straightforward. A possible
point of view suggested by the commutative situation is to describe
subspaces
by ideals and to glue algebras along ideals, using a pull-back (fibered product)
construction. This is the starting point of Budzy\'nski and Kondracki in
\cite{buko}, where coverings of $C^*$-algebras and locally trivial principal
fibre bundles over such covered algebras have been introduced.
In view of the fact that there are nontrivial algebras (e.g. the irrational
rotation
algebra) which have no nontrivial ideals, such an approach cannot reflect all
aspects of topological nontriviality of noncommutative algebras.
However, there are many algebras, in particular in the field of quantum
groups and quantum spaces, which have enough ideals for performing gluing
procedures, and it seems to be worthwhile to explore which kind of examples
may arise in this way.

If one wants to do differential geometry
in this scheme, the differential algebras (defining the differential
structure)
should also have a covering adapted to the covering of the underlying
algebra.
The construction of such adapted differential algebras will be one of the main
aims
of this article. It will be used in a subsequent paper, where we will consider
differential
structures  and connections on locally trivial quantum principal fibre
bundles
in the sense of \cite{buko}.

The present paper starts with the definition of coverings  of
algebras. Since we want to apply this notion also for differential algebras,
we cannot restrict ourselves to $C^*$-algebras, which leads to a difficulty being
absent there: It may happen that a natural gluing procedure fails to lead
from the collection of algebras corresponding to quantum subspaces defined
by a covering back to the original algebra. Coverings, which do not have
this
pathology, we call complete coverings. All the coverings of an algebra are
complete if the operations $+$ and $\bigcap$ between ideals are
distributive with respect to each other. If this general property of ideals
is not assumed, we can still give slightly weakened criteria for the
completeness of a covering. If an algebra is defined as a gluing, it has always 
a natural complete covering.

For differential algebras $\Gamma(B)=\bigoplus_{n\in\mathbbm{N}}
\Gamma^n(B)$ over an algebra $B$, we only consider
differential coverings being nontrivial with respect to $B$, i. e. coverings
consisting of differential ideals
whose components in degree zero form a nontrivial covering of $B$.
We show that for an algebra $B$ with given covering
$(J_i)_{i=1,
\ldots,n}$ and for given differential calculi $\Gamma(B_i)$ over the
algebras
$B_i=B/J_i$ corresponding to the ``local pieces'' of $B$ there exists a
unique
differential calculus $\Gamma(B)$ such that the natural projections
$\pi_i:B\lra B_i$ are differentiable and the kernels of the differential
extensions of the $\pi_i$ form a covering of $\Gamma(B)$. The covering 
completion with respect to this covering is in general only
locally a differential calculus. We also give a criterion assuring that the 
differential ideals generated by the ideals of a covering of an algebra $B$
form a covering of a given differential algebra $\Gamma(B)$.

The second part of the paper is devoted to an example. Gluing together
two copies of a quantum disc we obtain a $C^*$-algebra isomorphic to the
$C^*$-algebra of the Podle\'s sphere $S^2_{\mu c}$, $c>0$. This isomorphism,
already
mentioned in \cite{buko}, relies on the isomorphy of the disc algebra with
the $C^*$-algebra of the one-sided shift, and on a result of Sheu
\cite{sheu}
about the isomorphy of the gluing of two shift algebras by means of the
symbol
map and the Podle\'s sphere. We show that this $C^*$-algebra may also be
characterized as the $C^*$-closure of a ``polynomial'' algebra given in
terms
of generators and relations naturally arising from generators and relations
of the disc algebras via the gluing procedure. 
These generators should be considered as another set of ``coordinate
functions'' on the quantum sphere, which arise via a homeomorphism from
natural coordinates on a quantum version of a top of a cone. This is 
suggested by considering the spectra of the generators.

Finally, we construct, according to our general procedure, a differential
calculus on our ``quantum top'' out of two $U_{q^{1/2}}(sl_2)$-covariant
differential
calculi over the quantum discs. This differential calculus is also described in 
terms of relations between the generators and their differentials.
It is again a
$U_{q^{1/2}}(sl_2)$-covariant
differential calculus.\vspace{.5cm}

In the sequel, the word ``algebra'' always means an associative unital
algebra
over $\mathbbm{C}$. Ideals are always two-sided, and homomorphisms are homomorphisms
of algebras.

\section{ \label{cov}Coverings and gluings}

\indent

Let $B$ be an algebra and let $(J_i)_{i \in I}$ be a finite
family of ideals contained in B. Then the algebras $B_i$, $B_{ij}$ and $B_{ijk}$ are defined as the
factor algebras of $B$ with respect to the ideals $J_i$, $J_i+J_j$ and 
$J_i+J_j+J_k$.
The corresponding natural projections are denoted by
\begin{eqnarray*} \pi_i : B & \longrightarrow & B_i \\
\pi_{ij} : B & \longrightarrow & B_{ij} \\
\pi_{ijk} :B & \longrightarrow & B_{ijk}. \end{eqnarray*}
There are canonical surjective homomorphisms
\begin{eqnarray*} \pi^i_j :B_i & \longrightarrow & B_{ij} \\
\pi^i_{jk} :B_i & \longrightarrow & B_{ijk} \\
\pi^{ij}_k :B_{ij} & \longrightarrow & B_{ijk}.  \end{eqnarray*}
For example, $\pi^i_j(b+J_i)=b+J_i+J_j$.
Obviously, one can construct analogous factor algebras and surjective
homomorphisms
for a higher number of indices.
One easily shows $\pi^i_j \circ \pi_i=\pi^j_i \circ \pi_j=\pi_{ij}$ (and
similar formulas).
Furthermore, there are canonical
isomorphisms $B_{ij} \simeq B_i/\pi_i(J_j)$, which map $b+J_i+J_j
\in B_{ij}$ onto $ \pi_i(b) + \pi_i(J_j) \in B_i/\pi_i(J_j)$.
Note that the above definitions also mean $B_{ii}=B_i$, $\pi^i_i=\pi_i$,
etc..
\begin{de} Let $(J_i)_{i \in I}$ be
a finite family of ideals of an algebra $B$. $(J_i)_{i \in I}$ is
called covering of $B$ if \[ \bigcap_i J_i=\{0 \}. \]
A covering is called nontrivial if $J_i \not= \{0\}~~\forall i \in I$. 
\end{de}

For $C^*$-algebras and closed ideals, this definition was given by Budzy\'nski
and Kondracki \cite{buko}. For commutative $C^*$-algebras, this notion of
covering corresponds to coverings of the underlying topological space by
closed sets, the ideals just consisting of the functions vanishing on the
corresponding set.
We
want to use the definition also for differential algebras, which cannot be
made
$C^*$-algebras in an obvious way. Thus we are forced to stay in the general
algebraic context of our definition. As a consequence there may arise
difficulties with a reconstruction of the algebra from a covering by a
gluing
procedure which is always possible for $C^*$-algebras (see Proposition
\ref{cstar} below).
This is the motivation for introducing the notion of a complete covering.
\begin{de} \label{C} Let $B$ be an algebra and let $(J_i)_{i \in I}$ be a
covering of $B$. The algebra \begin{equation} \label{covcomp}
B_c:=\{ (a_i)_{i \in I} \in \bigoplus_{i \in I} B_i|
\pi^i_j(a_i)=\pi^j_i(a_j) \}. \end{equation}
is called the covering completion of $B$ with respect to $(J_i)_{i \in I}$.
The covering $(J_i)_{i \in I}$ is called complete if the injective
homomorphism \begin{equation} K: B \longrightarrow B_c  \end{equation}
defined by $K(a)= (\pi_i(a))_{i \in I}$ is surjective.
 $p_i :B_c \longrightarrow B_i$ denotes the restriction of the
canonical
projection $pr_i: \bigoplus_{j \in I} B_j \longrightarrow B_i$.
\end{de}

The name ``covering completion'' for $B_c$ is justified by the fact that
$(ker\,p_k)_{k \in I}$ is a complete covering of $B_c$, which will be
a special case of Proposition \ref{gluetogeth}.
$p_i$
is surjective, since $\pi_i$ is surjective and $\pi_i=p_i \circ K$.

Notice that the condition defining a complete covering is very similar to
one of the sheaf axioms:
It just says that a set of locally given objects which coincide ``on
intersections'' make up a global object. The other sheaf axiom, which says
that global objects, which coincide locally, also coincide globally,
corresponds to the injectivity of $K$, being true for any covering.

As shows the example below, there exist noncomplete coverings. On the other
hand, if an algebra has a covering, it also has a complete one:
If $(J_i)_{i\in I}$ is a nontrivial covering, consider ${\cal I}=
\{I'\subset I|\bigcap_{i\in I'}J_i=0\}$. Since the index set $I$ is 
finite, there exists $I'\in {\cal I}$ with minimal cardinality, i. e.
$card\,I'=min_{I''\in {\cal I}}card\,I''>1$ (since the covering is 
nontrivial). It follows that $\bigcap_{i\in I'}J_i=0$, $\bigcap_{i\in I''}
J_i\neq 0$, $\bigcap_{i\in I'\setminus I''}J_i\neq 0$ for any $I''\subset
I',~I''\neq I'$. Thus, $(\bigcap_{i\in I''}J_i,\bigcap_{i\in I'\setminus
I''}J_i)$ is a nontrivial covering.
Now, every two-element covering is complete:
\begin{pr} \label{K}
Let $B$ be an algebra and $J_1$, $J_2$ ideals of $B$. Then the
mapping $K:B\longrightarrow \{(a_1,a_2)\in B/J_1\bigoplus B/J_2|\pi^1_2(a_1)=
\pi^2_1(a_2)\}$ given by $K(b)=(\pi_1(b),\pi_2(b) )$ is surjective. In
particular, every covering consisting of two ideals is complete.
\end{pr}
Proof: One has to show that for every pair $(a_1,a_2) \in B_1 \bigoplus
B_2$ fulfilling $\pi^1_2(a_1)=\pi^2_1(a_2)$ there exists an element $a \in
B$
such that $\pi_1(a)=a_1$ and $\pi_2(a)=a_2$.
\\First we choose $\tilde{a_1},\;\tilde{a_2} \in B$ satisfying
$\pi_1(\tilde{a_1})=a_1$ and $\pi_2(\tilde{a_2})=a_2$. Clearly, there exist
elements
$r_1 \in J_1$ and $r_2 \in J_2$ such that \[ \tilde{a_1}=\tilde{a_2} + r_1 +
r_2
\] and one obtains the element $a$ by \[ a=\tilde{a_1}-r_1=\tilde{a_2} +
r_2. \] 
{} \hfill$\square$
\\\\

The following
proposition, which is a direct consequence of Theorems 17 and 18 (pages 279 and
280)
of \cite{zasa}, gives a sufficient condition for the completeness of
coverings:
\bpr
Assume that the operations $+$ and $\bigcap$ in the set of ideals of $B$ are
distributive with respect to each other, i. e. the set of ideals is a
distributive lattice with respect to these operations.
Then every covering of $B$ is complete.\label{zs}
\epr

This proposition is also true for subsets of the set of ideals of $B$ which
are closed under $+$ and $\cap$, with $+$ and $\cap$ distributive on the subset.

We will use similar arguments as in \cite{zasa} to prove criteria
for the completeness of a covering if the above condition is not assumed.
\begin{pr} \label{A} Let $B$ be an algebra and let $(J_i)_{i \in I}$ be a
covering of $B$. Assume that the index set is $I=\{1,2,...n\}$ and that the
ideals satisfy \[ \bigcap_{i=1,2...k-1}
(J_i+J_k)=(\bigcap_{i=1,2,...,k-1}J_i) +J_k,\;\; \forall k \in I. \]
Then $(J_i)_{i \in I}$ is complete.\end{pr}
Proof:  Notice that the condition $\pi^i_j(a_i)=
 \pi^j_i(a_j)$ just means $b_i-b_j\in J_i+J_j$ for $a_i=\pi_i(b_i)$,
 $a_j=\pi_j(b_j)$. Thus, in order to prove surjectivity of the map
 $a\longrightarrow (a+J_i)_i$, we have to show that from $b_i-b_j\in
J_i+J_j$
 for every pair of indices follows the existence of $b\in B$ with
 $b-b_i\in J_i$ (for all $i$). This is done inductively:
 Induction starts with Proposition \ref{K}. Assume now that we have
 found an $a\in B$ with $a-a_i\in J_i$ for $i=1,\ldots,k$. Then we have
 \[a-a_{k+1}=a-a_i+a_i-a_{k+1}\in J_i+J_i+J_{k+1}=J_i+J_{k+1},\]
 thus
 \[a-a_{k+1}\in
\bigcap\limits_{i=1}^k(J_i+J_{k+1})=(\bigcap\limits_{i=1}^kJ_i)+
 J_{k+1}.\]
 According to Proposition \ref{K} there exists $b\in B$ with $b-a\in\bigcap
 \limits_{i=1}^kJ_i$
 and $b-a_{k+1}\in J_{k+1}$. Therefore $b-a_i=b-a+a-a_i\in\bigcap
 \limits_{i=1}^kJ_i+J_i
 =J_i$, i. e. $b$ is the element we were looking for. \hfill $\square$
\begin{pr} \label{B} Let $B$ be an algebra and let $(J_i)_{i \in I}$ be a
complete covering of $B$. Then the family of ideals $(J_i)_{i \in I}$ has
the
property \[ \bigcap_{i \not= k}(J_i +J_k) =(\bigcap_{i \not= k}J_i) +
J_k\;\;\forall k
 \in I.\] \end{pr}
Proof: The inclusion $J_k+\bigcap_{i\neq k}J_i\subset\bigcap_{i\neq
k}(J_i+J_k)$
is true for subsets of a vector space. So we have to prove
$\bigcap_{i\neq k}(J_i+J_k)\subset J_k+\bigcap_{i\neq k}J_i$
for a complete covering. Completeness of the covering means that for
$(a_i)_{i\in I}\in\bigoplus_{i\in I}B_i$ with $\pi^i_j(a_i)=\pi^j_i(a_j)$
there exists a unique $a\in B$ with $\pi_i(a)=a_i$.
Let $I=\{1,...,n\}$ and assume $k=n$, without loss of generality.\\
Let $a\in\bigcap_{i<n}(J_i+J_n)$ and denote $a_i=\pi_i(a)$, $i=1,\ldots,n$.
Then we have $\pi^n_i(a_n)=\pi^n_i\circ\pi_n(a)=\pi_{in}(a)=0$ for $i<n$,
and
also $\pi^i_n(a_i)=\pi^i_n\circ\pi_i(a)=\pi_{ni}(a)=0$ for $i<n$,
from which we can conclude $(0,\ldots,a_n)\in B_c$ and
$(a_1,\ldots,a_{n-1},0)\in B_c$. Obviously,
$(a_1,\ldots,a_{n-1},0)\in ker\, p_n$ and $(0,\ldots,0,a_n)\in
\bigcap_{i<n}ker p_i$.
On the other hand, $K:B\longrightarrow B_c$ is by assumption an algebra
isomorphism
and one easily verifies that it maps $J_i$ onto $ker\, p_i$. Therefore,
there
are $b\in J_n$, $c\in \bigcap_{i<n}J_i$ such that
$K(b)=(a_1,\ldots,a_{n-1},0)$
and $K(c)=(0,\ldots,0,a_n)$, and we have \[a=K^{-1}(K(b)+K(c))=b+c\in
J_n+\bigcap_{i<n}J_i.\]
{}\hfill$\square$
\begin{pr}\label{3ue}
For a covering $(J_1,J_2,J_3)$ consisting of three ideals the following
conditions are equivalent:\\
(i) $(J_1,J_2,J_3)$ is a complete covering.\\
(ii) $(J_i+J_k)\bigcap(J_j+J_k)=J_i\bigcap J_j + J_k$ for one permutation
$i,j,k$
of $1,2,3$.\\
(iii) $(J_i+J_k)\bigcap(J_j+J_k)=J_i\bigcap J_j + J_k$ for every permutation
$i,j,k$
of $1,2,3$.
\end{pr}
Proof: Is an obvious combination of the two foregoing propositions according
to the
scheme (i) $\Rightarrow$ (iii) $\Rightarrow$ (ii) $\Rightarrow$ (i). The
right
$\Rightarrow$ is possible only for a three-element covering, because only
then the
conditions of Proposition \ref{A} reduce to exactly one of the conditions of
Proposition \ref{B}. \hfill$\square$
\\\\

Let us also note that a covering is always complete if the ideals are
coprime,
i. e. if $J_i+J_j=B$, $i\neq j$ (\cite{atmac}). However, in this case there
is no gluing at all.
\\\\Example of a noncomplete covering:
\\Consider the algebra
\[B=\mathbbm{C}<x,y,z>/J,\]
where $J$ is the ideal generated by the elements $xy,yx,xz,zx,yz,zy$.
It follows that
\[\{1,x,x^2,\ldots,y,y^2,\ldots,z,z^2,\ldots\}\]
is a linear basis of $B$. Consider the ideals $J_1$, $J_2$, $J_3$
generated by $x-y$, $x-z$, $y-z$, respectively. Obviously,
\[J_1=\mathbbm{C}(x-y)+\mathbbm{C} x^2+\mathbbm{C} x^3+\ldots +\mathbbm{C} y^2+\mathbbm{C} y^3+\ldots,\]
\[J_2=\mathbbm{C}(x-z)+\mathbbm{C} x^2+\mathbbm{C} x^3+\ldots +\mathbbm{C} z^2+\mathbbm{C} z^3+\ldots,\]
\[J_3=\mathbbm{C}(y-z)+\mathbbm{C} y^2+\mathbbm{C} y^3+\ldots+\mathbbm{C} z^2+\mathbbm{C} z^3+\ldots.\]
Moreover,
\[\mathbbm{C}(x-y)+\mathbbm{C}(x-z)=\mathbbm{C}(x-y)+\mathbbm{C}(y-z)=\mathbbm{C}(x-z)+\mathbbm{C}(y-z),\]
whereas
\[\mathbbm{C}(x-y)\cap\mathbbm{C}(x-z)=\mathbbm{C}(x-y)\cap\mathbbm{C}(y-z)=\mathbbm{C}(x-z)\cap\mathbbm{C}(y-z)
=
\{0\}.\]
Therefore, we obtain
\[J_1\cap J_2\cap J_3=\{0\},\]
i.e. $(J_1,J_2,J_3)$ is a covering of $B$. For the sums and intersections
of
two of the three ideals we get
\[J_1+J_2=\mathbbm{C}(x-y)+\mathbbm{C}(x-z)+\mathbbm{C} x^2+\ldots+\mathbbm{C} y^2+\ldots+\mathbbm{C}
z^2+\ldots,\]
\[J_1+J_3=J_2+J_3=J_1+J_2,\]
\[J_1\cap J_2=\mathbbm{C} x^2+\mathbbm{C} x^3+\ldots,\]
\[J_1\cap J_3=\mathbbm{C} y^2+\mathbbm{C} y^3+\ldots,\]
\[J_2\cap J_3=\mathbbm{C} z^2+\mathbbm{C} z^3+\ldots.\]
We conclude that
\[(J_1+J_3)\cap(J_2+J_3)=J_1+J_3,\]
whereas
\[(J_1\cap J_2)+J_3=\mathbbm{C}(y-z)+\mathbbm{C} x^2+\ldots+\mathbbm{C} y^2+\ldots+\mathbbm{C} z^2+
\ldots,\]
which is strongly contained in $J_1+J_3$. Similarly,
all other possible equalities of
Propositions \ref{A} and \ref{B} are not
 satisfied, which means that the covering is not complete.

Notice that we could have introduced the additional relations $x^n=0$,
$y^n=0$,
$z^n=0$, $n\geq 2$, for example, without changing the situation essentially.
For $n=2$ we arrive at a pure vector space situation (three subspaces such
that
the sum of any two contains the third). Admittedly, in this case the
covering
is reducible, already two of the three ideals form a
covering.\\[.3cm]
Second example of a noncomplete covering:\\
Let ${\Bbb C}[x,y]$ be the algebra of polynomials in the (commuting) 
indeterminates $x$ and $y$. Consider the principal ideals $J_1=(x)$, $J_2=(y)$,
$J_3=(x-y)$. One easily verifies that
\[J_1+J_2=J_1+J_3=(J_1+J_2)\cap(J_1+J_3)=\mbox{polynomials without constant term},\]
\[J_1+J_2\cap J_3=(x)+(y(x-y)),\]
i. e. $(J_1+J_2)\cap(J_1+J_3)\neq J_1+J_2\cap J_3$. Moreover,
\[J_1\cap J_2\cap J_3=(xy(x-y)).\]
Therefore, going to the factor algebra $A={\Bbb C}[x,y]/J$, where $J$ is the
ideal generated by monomials of at least third degree, we obtain a noncomplete 
covering $(J_1,J_2,J_3)$. The triple $(x-y+J_1,x-y+J_2,x+J_3)$ is an 
element of $A_c$ which has no preimage in $A$.
%
\begin{pr} \label{cstar} Any covering of a $C^*$-algebra consisting of
closed
ideals is complete.
\end{pr}
Proof: The closed ideals of a $C^*$-algebra form a distributive lattice with
respect to the operations $+$ and $\cap$,
which in turn follows from the fact that in this case the product of closed
ideals coincides with their intersection, see \cite{di1}, 1.9.12.a..
\hfill$\square$
\\\\

The proof of the following proposition can be found in \cite{buko}.
\begin{pr} A $C^*$-algebra which admits a faithful irreducible
representation does not admit any nontrivial covering consisting of closed
ideals.\end{pr}

In particular, the algebra $B(H)$ of bounded operators on a Hilbert space
$H$
 does not admit a nontrivial covering. The same is obviously true for any
 simple algebra.
\\\\

A general method to construct algebras possessing a complete covering
is given by a gluing procedure:
\bde\label{glu}
Assume that there are given finite
families of algebras $(B_i)_{i \in I}$ and $(B_{ij})_{i,j \in I}$
and surjective homomorphisms $\pi^i_j: B_i \longrightarrow B_{ij}$,
where $B_{ij}=B_{ji}$, $B_{ii}=B_i$ and  $\pi^i_i=id$.
Then the algebra
\[ \oplus_{\pi^i_j}B_i:= \{ (a_i)_{i \in I} \in \bigoplus_{i \in I}
B_i|~\pi^i_j(a_i)=\pi^j_i(a_j) \}\]
is called gluing of the algebras $B_i$ with respect to the $\pi^i_j$.
\ede

For $I=\{1,2\}$, this is known as a pull-back or a fibered product of $B_1$
and $B_2$. The covering completion $B_c$ of an algebra $B$ with covering 
$(J_i)_{i\in I}$ is just the gluing of the $B/J_i$ with respect to the natural
maps $\pi^i_j$.
We will show that $(ker\, p_i)_{i \in I}$, where
$p_i: \oplus_{\pi^i_j}B_i \longrightarrow B_i$ are the restrictions of the canonical
projections,
 is a complete
covering of $B$. However, the $p_i$ are in
general not surjective.
In the classical situation, this would mean that the sets, which are glued
together, are not embedded in the global object.
\begin{lem} \label{phiijhom}
 Let $B=\oplus_{\pi^i_j}B_i$ and $A=\oplus_{\eta^i_j}A_i$ be gluings as in
Definition \ref{glu}. Moreover, let\\
$\phi_i:A_i\longrightarrow B_i$ be algebra homomorphisms. \\
Assume that there exist algebra homomorphisms $\phi_{ij}=\phi_{ji}:A_{ij}
\longrightarrow B_{ij}$ such that
\begin{equation}
\phi_{ij}\circ\eta^i_j=\pi^i_j\circ\phi_i, \; i\neq j.\label{phij}
\end{equation}
Then we have $(\oplus_i\phi_i)(A)\subset B$ .
\end{lem}
Proof:\\
$(\oplus_i\phi_i)(A)\subset B$ means
$\eta^i_j(a_i)=\eta^j_i(a_j)\;\Rightarrow\;\pi^i_j(\phi_i(a_i))=\pi^j_i(\phi
_j
(a_j))$. Assuming the existence of $\phi_{ij}$ with (\ref{phij}) it follows
from
$\eta^i_j(a_i)=\eta^j_i(a_j)$ that
$\pi^i_j\circ\phi_i(a_i)=\phi_{ij}\circ\eta^i_j
(a_i)=\phi_{ij}\circ\eta^j_i(a_j)=\pi^j_i\circ\phi_j(a_j)$.
\hfill$\square$\\[.1cm]
\\

Notice that the $\phi_{ij}$ fulfilling (\ref{phij}) exist if and only if
$ker(\eta^i_j)\subset ker(\pi^i_j\circ
\phi_i),\;i\neq j$. Moreover, for the lemma it is not necessary that
$\pi^i_j$ and $\eta^i_j$ are surjective.
\begin{pr} \label{gluetogeth}
Let $B=\oplus_{\pi^i_j}B_i$.
\\
Then $(ker\,p_i)_{i \in I}$ is a complete covering of B. 
\end{pr}
Proof:	It is obvious that $(ker\,p_i)_{i\in I}$ is a covering of $B$. The
covering
completion $B_c$ with respect to this covering is defined as
\[B_c=\{(b_i)_{i\in I}\in\bigoplus_i
B/ker\,p_i|~\eta^i_j(b_i)=\eta^j_i(b_j)\},\]
where $\eta^i_j:B/ker\,p_i \longrightarrow B/(ker \,p_i+ker \,p_j)$ are
the
canonical maps $(b_k)_{k\in I}+ker \,p_i \longrightarrow (b_k)_{k\in I}
+ker\,p_i+ker\,p_j$.

Let $\phi_i:B/ker\,p_i \longrightarrow B_i$ be defined as
$(b_k)_{k\in I}+ker \,p_i \longrightarrow b_i$. Obviously, the $\phi_i$ are
well defined and injective.\\
It is now sufficient to show that $(\oplus_i\phi_i)(B_c)\subset B$ and
$K\circ
(\oplus_i\phi_i)=id_{B_c}$, where $K:B\longrightarrow B_c$ is the canonical
embedding $(b_j)_{j\in I}\longrightarrow ((b_j)_{j\in I}+ker\,p_i)_{i\in
I}$. The latter is a trivial verification. In order to
verify the first claim, define
$\phi_{ij}:B/(ker \,p_i+ker \,p_j)\longrightarrow B_{ij}$ by
$\phi_{ij}\circ\eta^i_j=\pi^i_j\circ\phi_i$, i. e.
$\phi_{ij}((b_k)_{k\in I}+ker\,p_i+ker \,p_j):=\pi^i_j(b_i)$.
According to the remark after Lemma \ref{phiijhom}, $\phi_{ij}$ is well defined:
$ker(\eta^i_j)=\{(b_k)_{k\in I}+ker \,p_i\in B/ker \,p_i|~b_j=0\}\subset
ker\pi^i_j\circ\phi_i=\{(b_k)_{k\in I}+ker\,p_i \in B/ker \,p_i|~b_i\in
ker\pi^i_j\}$, since from $b_j=0$ follows $\pi^i_j(b_i)=\pi^j_i(b_j)=0$.
Thus, Lemma \ref{phiijhom} proves the claim.\hfill$\square$\\[.1cm]
\\

Notice that the implication $b_i=0\Rightarrow\pi^j_i(b_j)=0
$ also means $p_j(ker\,p_i)\subset ker(\pi^j_i)$.\\[.3cm]

Possible nonsurjectivity of $p_i$ is reflected in nonsurjectivity of the
$\phi_i$ appearing in the foregoing proof. In the classical situation of
algebras of functions over compact spaces, this would mean that the space underlying the algebra $B_i$ is not injectively mapped onto the space
underlying $B/ker\,p_i$. With other words, the $\phi_i$ would encode a 
gluing
of $B_i$ with itself. The gluing of the $B/ker\,p_i$ does not lead to a
further self-gluing, in contrast to the gluing of the $B_i$, where all
the gluing is done ``in one step''.\\[.3cm]

Proposition \ref{gluetogeth} also has the consequence that the covering
completion $B_c$ of an algebra $B$ with covering $(J_i)_{i\in I}$ has the
complete covering $(ker\,p_i)_{i\in I}$.\\

If $B_i$ and $B_{ij}$ are $C^*$-algebras, the
kernels of the homomorphisms $\pi^i_j$ are closed ideals. The distributivity
of $+$ and $\bigcap$ on the set of closed ideals leads to the following
sufficient condition for the surjectivity of the projections $p_i$:
\begin{pr} Let the algebras $B_i$ and $B_{ij}$ be $C^*$-algebras. Assume
that
the homomorphisms $\pi^i_j$ have the following properties:
\\\\1. \begin{equation} \label{pishit1} \pi^i_j(ker \pi^i_k)=
\pi^j_i(ker \pi^j_k);~~\forall i,j,k \in I \end{equation}
2. The isomorphisms $\pi^{ij}_k: B_i/(ker \pi^i_j + ker \pi^i_k)
\longrightarrow
B_{ij}/ \pi^i_j(ker \pi^i_k)$ defined by
\[ \pi^{ij}_k(f + \ker \pi^i_j + ker\pi^i_k)=\pi^i_j(f) + \pi^i_j(ker
\pi^i_k)
\] fulfill \begin{equation} \label{pishit} {\pi^{ik}_j}^{-1} \circ
\pi^{ki}_j={\pi^{ij}_k}^{-1} \circ \pi^{ji}_k \circ {\pi^{jk}_i}^{-1} \circ
\pi^{kj}_i. \end{equation}
Then the homomorphisms $p_i$ are surjective. 
\end{pr}
Remark 1: In the classical situation, i. e. $B_i=C(U_i)$, $B_{ij}=C(U_{ij})$,
$U_i,~U_{ij}$ compact spaces, $\pi^i_j$ is the pull-back of an embedding
$\iota^i_j:U_{ij}\lra U_i$, and $ker\pi^i_j$ are the functions vanishing on
$\iota^i_j(U_{ij})\subset U_i$. Condition (\ref{pishit1}) says that
$\iota^i_j(U_{ij})\bigcap \iota^i_k(U_{ik})$ is homeomorphic to
$\iota^j_i(U_{ij})\bigcap \iota^j_k(U_{jk})$ for every triple $i,j,k$.
$\pi^{ij}_k$ is the pull-back of the restriction of $\iota^i_j$ to
${\iota^i_j}^{-1}(\iota^i_k(U_{ik}))$. (\ref{pishit}) is the natural
compatibility
condition for these restrictions.\\[.3cm]
Remark 2: If $B$ is an algebra with covering $(J_i)_{i\in I}$, then the 
homomorphisms $\pi^i_j:B_i=B/J_i\lra B_{ij}=B/(J_i+J_j)$ defining the
covering completion $B_c$ satisfy the assumptions of the proposition.\\[.3cm]
Proof: Assume that homomorphisms $\pi^i_j$ satisfying
(\ref{pishit1})
and (\ref{pishit}) are given. Let the isomorphisms $\phi^k_{ij}: B_j/(ker \pi^j_i +
ker
\pi^j_k) \longrightarrow B_i/(ker \pi^i_j + ker \pi^i_k)$ be defined by
$\phi^k_{ij}:= {\pi^{ij}_k}^{-1} \circ {\pi^{ji}_k}$.
The $\phi^k_{ij}$ satisfy (see formula (\ref{pishit}) )
$\phi^j_{ik}=\phi^k_{ij} \circ \phi^i_{jk}$ and
${\phi^k_{ij}}^{-1}=\phi^k_{ji}$.
To prove the surjectivity of the projection $p_i$ one has
to show that for all $f \in B_i$ there exists a family $(f_k)_{k \in I} \in
B$
such that $p_i((f_k)_{k \in I})=f$. Suppose that the index set is
$I=\{1,2,3,...,n\}$.
It is sufficient to consider the case $i=1$.
For $f \in B_1$ there exist $f_2 \in B_2$ such that
$\pi^1_2(f)=\pi^2_1(f_2)$ and $f_3 \in B_3$
such
that $\pi^1_3(f)=\pi^3_1(f_3)$. It follows that \begin{eqnarray*}
f + ker \pi^1_2 + ker \pi^1_3 &=& \phi^3_{12}(f_2 + ker\pi^2_1 + ker
\pi^2_3) \\ &=& \phi^2_{13}(f_3 + ker \pi^3_1 + ker \pi^3_2)
\end{eqnarray*}
and, with $\phi^1_{32}={\phi^2_{13}}^{-1}\ci\phi^3_{12}$,
 \[ f_3 + ker \pi^3_1 + ker \pi^3_2= \phi^1_{32}(f_2 +ker\pi^2_1 + ker
\pi^2_3), \] thus $\pi^3_2(f_3)-\pi^2_3(f_2)=r_{23} \in \pi^3_2(ker
\pi^3_1)$.
Choose $\tilde{r}_{23} \in ker \pi^3_1$ such that
$\pi^3_2(\tilde{r}_{23})=r_{23}$. Then $f^1_3 := f_3 - \tilde{r}_{23}$
satisfies $\pi^1_3(f)= \pi^3_1(f^1_3)$ and $\pi^2_3(f_2)= \pi^3_2(f^1_3)$.
\\Now assume that a family 
$(f_j)_{j=1,...,k},~~f_j \in
B_j$ with $f_1=f$
and $\pi^i_j(f_i)= \pi^j_i(f_j),~~\forall i,j=1,...,k$ has been found.
 Then there exists $f^1_{k+1} \in
B_{k+1}$ with $\pi^1_{k+1}(f)= \pi^{k+1}_1(f^1_{k+1})$.
 Assume now that for fixed $i \in \{1,...,k-1\}$ there is given
$f^i_{k+1} \in B_{k+1}$ which satisfies $\pi^j_{k+1}(f_j)=
\pi^{k+1}_j(f^i_{k+1})$ for $j=1,...,i$.
It follows that there exists
$f^{i+1}_{k+1} \in B_{k+1}$ which satisfies  $\pi^j_{k+1}(f_j) =
\pi^{k+1}_j(f^{i+1}_{k+1}),~~
\forall j=1,...,i+1$: There are the identities
\begin{eqnarray*}
f_j + ker \pi^j_{k+1} + ker \pi^j_{i+1}&=&\phi^{i+1}_{j,k+1}(f^i_{k+1} +ker
\pi^{k+1}_j + ker \pi^{k+1}_{i+1}),~ \forall j=1,...,i, \\
&=& \phi^{k+1}_{j,i+1}(f_{i+1} + ker
\pi^{i+1}_j + ker \pi^{i+1}_{k+1}),~~ \forall j=1,...,i, \end{eqnarray*}
which lead to
\[ f^i_{k+1} + ker \pi^{k+1}_j + ker \pi^{k+1}_{i+1}=
\phi^j_{k+1,i+1}(f_{i+1} +
ker \pi^{i+1}_j + ker \pi^{i+1}_{k+1}),~~ \forall j=1,...,i, \] and it
follows that
\[ \pi^{k+1}_{i+1}(f^i_{k+1}) - \pi^{i+1}_{k+1}(f_{i+1}) = r_{i+1,k+1} \in
\bigcap_{j=1,...,i} \pi^{k+1}_{i+1} (ker \pi^{k+1}_j). \]
Because of
\[ \bigcap_{j=1,...,i} (ker \pi^{k+1}_j + ker \pi^{k+1}_{i+1})=
(\bigcap_{j=1,...,i} ker \pi^{k+1}_j) + ker \pi^{k+1}_{i+1}, \]
in the case of $C^*$- algebras, applying $\pi^{k+1}_{i+1}$ one obtains
$\bigcap_{j=1,...,i} \pi^{k+1}_{i+1} (ker \pi^{k+1}_j)\\
=\pi^{k+1}_{i+1} (\bigcap_{j=1,...,i} ker \pi^{k+1}_j) $. Thus one finds
$\tilde{r}_{i+1,k+1} \in \bigcap_{j=1,...,i} ker
\pi^{k+1}_j$, such that $\pi^{k+1}_{i+1}(\tilde{r}_{i+1,k+1})\newline
=r_{i+1,k+1}$,
and $f^{i+1}_{k+1}=f^i_{k+1}-\tilde{r}_{i+1,k+1}$
satisfies $\pi^{k+1}_{j}(f^{i+1}_{k+1})=\pi^j_{k+1}(f_j),~~\forall
j=1,...,i+1$.
\\This means that there exists $f_{k+1} \in B_{k+1}$ satisfying
\[\pi^{k+1}_j(f_{k+1})=\pi^j_{k+1}(f_j),~~\forall j=1,...,k.\] 
Continuing this procedure one obtains a family $(f_i)_{i \in I} \in B$ with
\\$p_1((f_i)_{i \in I})=f$. Thus $p_1$ is surjective. \hfill$\square$
\\\\ If only two algebras are glued together
the
projections $p_1$ and $p_2$ are always surjective.
\section{\label{adapt}Adapted differential structures on algebras with
covering}
\begin{de} A differential algebra $\Gamma (B)$ over an
algebra $B$ is an $\mathbbm{N}-$graded algebra, i.e.
\begin{eqnarray*}
 \Gamma (B)&=&\bigoplus_{n \in \mathbbm{N}} \Gamma^n(B) \\
 \Gamma^n(B) \Gamma^m(B) &\subset & \Gamma^{m+n}(B), 
\end{eqnarray*}
with
 \[ \Gamma^0(B)=B, \]
 which is equipped with a differential, i.e. a
linear
 map $d$ of
$\Gamma(B)$ fulfilling
\begin{eqnarray*}
   d(\Gamma^n(B)) &\subset & \Gamma^{n+1}(B)  \\
 d(\rho  \eta )&=&(d \rho ) \eta + (-1)^n \rho  d\eta,
 \;\; \rho \in \Gamma^n(B),\; \eta \in \Gamma(B)\\
d^2&=&0. 
\end{eqnarray*}
A differential ideal $J \subset \Gamma(B)$ is an ideal of the algebra
$\Gamma(B)$
such that \begin{eqnarray} pr_i(J) &\subset& J \\ dJ & \subset& J,
\end{eqnarray}
where $pr_i: \Gamma(B) \longrightarrow \Gamma^i(B)$ is the canonical
projection.\\
A homomorphism $\phi:\Gamma(B)\lra\Gamma(A)$ of differential algebras
(of degree $0$) is an algebra homomorphism with
$\phi\ci d=d\ci\phi$ and $pr_i\ci\phi=\phi\ci pr_i$.
\end{de}

Differential ideals are in bijective correspondence with kernels of
surjective homomorphisms of differential algebras.
\begin{de} 
A differential algebra $ \Gamma(B)$ over an algebra
$B$ is called differential calculus, if every element $\rho \in \Gamma^n(B)$
has
the general form \begin{equation}
\rho =\sum_k a^k_0da^k_1 ... da^k_n, \;\;a^k_i \in B, \end{equation} i. e.
if
$B$ and $dB$ generate $\Gamma(B)$ as an algebra.
 \end{de}

If $\Gamma(B)$ is a differential algebra and $J \subset \Gamma(B)$ is a
differential ideal, $\Gamma(B)/J$ is a differential algebra over
$B/pr_0(J)$.
If $\Gamma(B)$ is a differential calculus, $\Gamma(B)/J$ is a differential
calculus
over $B/pr_0(J)$.

For any algebra there exists the universal differential calculus 
$\Omega(B)$ over $B$.
As is
well known, every differential calculus $\Gamma(B)$ over $B$ corresponds to
a
unique differential ideal $J(B)\subset\Omega(B)$ such that $\Gamma(B)\simeq
\Omega(B)/J(B)$.

For the next definition see also \cite{pfscha}.
\begin{de} Let $\Psi:A\longrightarrow B$ be a homomorphism of algebras, and
let $\Gamma
(A)$, $\Gamma(B)$ be differential algebras over $A, B$. $\Psi$ is called
differentiable with respect to $\Gamma(A)$ and $\Gamma(B)$ if there exists
a homomorphism of graded algebras
$\Psi_\Gamma:\Gamma(A)\longrightarrow\Gamma(B)$
such that 
\[\Psi_\Gamma\circ d=d\circ\Psi_\Gamma,\]
\[\Psi_\Gamma|_A=\Psi,\]
i. e. if there exists an extension of $\Psi$ to a homomorphism of
differential algebras. $\Psi_\Gamma$ is called extension of $\Psi$
with respect to $\Gamma(A)$ and $\Gamma(B)$.
\end{de}

The following statements are well known or obvious:

If $\Gamma(A)$ is a differential calculus, the extension $\Psi_\Gamma$ is
uniquely determined by
\begin{equation} 
\Psi_{\Gamma}(a_0da_1 ...
da_n)=\Psi(a_0)d\Psi(a_1) ...  d\Psi(a_n). 
\end{equation}
If $\Gamma(A)=\Omega(A)$, the universal differential calculus over $A$,
$\Psi_\Gamma$ always exists (see \cite{sei}). We denote it sometimes by
$\Psi_{\Omega\rightarrow\Gamma}$. If $\Psi$ is surjective,
$\Psi_{\Omega\rightarrow\Gamma}$ is surjective if and only if $\Gamma(B)$
is a differential calculus. The extension of $\Psi$ with respect to the
universal
differential calculi $\Omega(A)$ and $\Omega(B)$ will be denoted by
$\Psi_{\Omega}$. If both $\Gamma(A)$ and $\Gamma(B)$ are differential
calculi, $\Psi$ is differentiable with respect to
$\Gamma(A)\simeq\Omega(A)/
J(A)$ and $\Gamma(B)\simeq \Omega(B)/J(B)$ if and only if
$\Psi_\Omega(J(A))\subset J(B)$.
%
\begin{de} Let $\Gamma(B)$ be a differential algebra.
A covering $(J_i)_{i\in I}$ of $\Gamma(B)$
is called differentiable if the $J_i$ are differential ideals.\\
A differentiable covering is complete if it is complete as a covering.
\end{de}
\begin{pr} 
Let $(J_i)_{i \in I}$ be a differentiable covering of the
differential algebra $\Gamma(B)$. Then $(pr_0(J_i))_{i \in I}$ is a covering
of $B$. 
\end{pr}
Proof:  
\[ \bigcap_i pr_0(J_i) =\{ a \in B| a\in \bigcap_i J_i \}=\{0\}. \] 
\hfill$\square$
\begin{de} If the covering of $B$ induced by a differentiable covering of
$\Gamma(B)$ is nontrivial the differential covering is said to be nontrivial
with respect to $B$. 
\end{de}

One easily finds nontrivial differential coverings of the algebra of usual
differential forms on a manifold, which are trivial with respect
to the algebra of smooth functions on the manifold, e. g. just
putting the zeroth degree to 0. The above definition is used to avoid such
cases.

\begin{de} A differential algebra $\Gamma(B)$ with a complete differentiable
covering $(J_i)_{i \in I}$ which is nontrivial with respect to $B$ is called
LC differential algebra, if the factor algebras $\Gamma(B)/J_i$ are
differential calculi over the algebras $B/pr_0(J_i)$. 
\end{de}

LC differential algebras will naturally arise from differential structures
on
locally trivial quantum principal fibre bundles (see \cite{coma2}).

\begin{de} Let $(B,(J_i)_{i\in I})$ be an algebra with covering,
let $B_i=B/J_i$, let $\pi_i:B\lra B_i$ be the natural surjections, and
let $\Gamma(B)$
and $\Gamma(B_i)$ be differential calculi such that $\pi_i$ are 
differentiable and $(ker{\pi_i}_\Gamma)_{i\in I}$ is a covering
of $\Gamma(B)$. Then $(\Gamma(B),(\Gamma(B_i))_{i\in I})$ is called adapted
to
$(B,(J_i)_{i\in I})$.
\end{de}

If $(J_i)_{i\in I}$ is a nontrivial covering in this situation, $(ker\,
{\pi_i}_\Gamma)_{i\in I}$ is nontrivial with respect to $B$, since
$pr_0(ker\,{\pi_i}_\Gamma)=J_i$.

In the classical case, where the $\pi_i$ are the pull-backs of embeddings of
closed submanifolds $M_i$ into a manifold $M$ and the $\Gamma$'s are usual
differential forms, the $\pi_i$ are obviously differentiable, the
${\pi_i}_\Gamma$ being the pull-backs on forms, and $ker{\pi_i}_\Gamma$ are
the
differential forms vanishing on $M_i$, thus obviously forming a covering.
\begin{pr} \label{diff} 
Let $(B,(J_i)_{i \in I})$ be an algebra with
covering, and let $\Gamma(B_i)$ be differential calculi over the
algebras $B_i$. Up to isomorphy there exists a unique differential
calculus
$\Gamma(B)$ such that $(\Gamma(B),(\Gamma(B_i))_{i \in I}))$ is adapted to
$(B,(J_i)_{i \in I})$. 
\end{pr}
Proof: As noted above there exist the extensions
$\pi_{i_{\Omega \rightarrow \Gamma}} :\Omega(B) \longrightarrow \Gamma(B_i)$
of
$\pi_i$ defined by \[ \pi_{i_{\Omega \rightarrow \Gamma}}
(a_0da_1...da_n):=\pi_i(a_0)d\pi_i(a_1) ...d \pi_i(a_n). \]
Then $J(B):= \bigcap_i ker \pi_{i_{\Omega \rightarrow \Gamma}}$ is a
differential ideal in $\Omega(B)$. Because of
$pr_0(J(B))=J(B) \bigcap B =0$, $\Gamma(B):=\Omega(B) / J(B)$ is a differential
calculus over $B$. The extensions $\pi_{i_{\Gamma}}$ of $\pi_i$ exist
and the pair $(\Gamma(B),( \Gamma(B_i))_{i \in I})$
is adapted to $(B,(J_i)_{i \in I})$.
\\Let $\tilde{\Gamma}(B)=\Omega(B)/\tilde{J}$ be another differential
calculus such
that $(\tilde{\Gamma}(B),(\Gamma(B_i))_{i\in I})$ is adapted to
$(B,(J_i)_{i\in
I})$. Let $\pi_{\Omega,
\tilde{\Gamma}}:\Omega(B)\longrightarrow\tilde{\Gamma}(B)$
denote the canonical quotient map. Then differentiability of $\pi_i$ with
respect
to $\tilde{\Gamma}(B)$ and $\Gamma(B_i)$ means that there exist
$\pi_{i_{\tilde{\Gamma}}}:\tilde{\Gamma}(B)\longrightarrow\Gamma(B_i)$ such
that
$\pi_{i_{\tilde{\Gamma}}}\circ\pi_{\Omega,\tilde{\Gamma}}=\pi_{i_{\Omega
\rightarrow
\Gamma}}$. Therefore we have
$\tilde{J}=ker\pi_{\Omega,\tilde{\Gamma}}\subset
ker\pi_{i_{\Omega\rightarrow\Gamma}}$ $\forall i$, i. e.
$\tilde{J}\subset\bigcap_{i}
ker\pi_{i_{\Omega\rightarrow\Gamma}}=J$.
On the other hand, if $\gamma\in J\setminus\tilde{J}$, then
$\gamma+\tilde{J}$
is a nonzero element of $\bigcap_i ker\pi_{i_{\tilde{\Gamma}}}$, since
$\pi_{i_{\tilde{\Gamma}}}(\gamma+\tilde{J})=\pi_{i_{\Omega\rightarrow\Gamma}
}(\gamma)
=0$ $\forall i$.
\hfill$\square$
\\

For a given differential calculus
$\Gamma(B)=\Omega(B)/J(B)$, the induced differential calculi
$\Gamma(B_i)=\Omega(B_i)/\pi_{i_{\Omega}}(J(B))$ in general do not form a
pair
$(\Gamma(B),(\Gamma(B_i))_{i \in I})$ which is adapted to $(B,(J_i)_{i
\in I})$:
\begin{pr} 
Let $\Gamma(B):= \Omega(B) /J(B)$ be a differential calculus over
the
algebra $B$, and let $(J_i)_{i \in I}$ be a covering of $B$.
\\The pair $(\Gamma(B),(\Omega(B_i) /\pi_{i_{\Omega}}(J(B))))_{i\in I})$ is
adapted to $(B,(J_i)_{i \in I})$ if and only if the differential ideal
$J(B)$ has the property
\begin{equation} \label{adapti} J(B)=\bigcap_{i \in I}(J(B) + ker
\pi_{i_{\Omega}}). \end{equation} 
\end{pr}
Remark 1: Obviously, not every differential ideal has this property. For
example, in the case of universal differential calculi,
where $J(B)=0$, the differential ideal
$\bigcap_i ker \pi_{i_{\Omega}}$ contains elements of
$J_{\sigma(1)}J_{\sigma(2)}...J_{\sigma(n-1)}dJ_{\sigma(n)}$ in the first
degree,
where $\sigma$ is any permutation, thus condition (\ref{adapti}) is not
satisfied.\\\\
Remark 2: $\Gamma(B_i)=\Omega(B_i)/\pi_{i_{\Omega}}(J(B))$ is canonically
isomorphic to $\Gamma(B)/(J_i,dJ_i)$, where $(J_i,dJ_i)$ is the differential
ideal generated by $J_i$ in $\Gamma(B)$. $((J_i,dJ_i))_{i\in I}$ is a 
covering of $\Gamma(B)$ iff (\ref{adapti}) is satisfied.
\\\\Proof: First consider a differential ideal $J(B)$ fulfilling
condition
(\ref{adapti}), and differential calculi $\Gamma(B_i):=\Omega(B_i)/
\pi_{i_{\Omega}}(J(B))$. One has to prove $\bigcap_{i \in I} ker
\pi_{i_{\Omega \rightarrow \Gamma}}=J(B)$.
\\Since $ker \pi_{i_{\Omega \rightarrow
\Gamma}}=\pi^{-1}_{i_{\Omega}} (
\pi_{i_{\Omega}}(J(B)))=J(B) + ker \pi_{i_{\Omega}}$, and since $J(B)$
fulfills
(\ref{adapti}), one direction of the assertion is proved.
\\[.2cm]
Now let $(\Gamma(B),(\Gamma(B_i))_{i \in I})$ be adapted, i.e. $J(B)=
\bigcap_{i \in I} ker \pi_{i_{\Omega \rightarrow \Gamma}}$.
$J(B)\subset \bigcap_{i\in I}(J(B)+ker\,\pi_{i_\Omega})$ is trivially 
true.
Conversely, $ \bigcap_{k \in I} ker \pi_{k_{\Omega \rightarrow
\Gamma}}
\subset ker \pi_{i_{\Omega \rightarrow \Gamma}}$ and $ker \pi_{i_{\Omega}}
\subset ker \pi_{i_{\Omega \rightarrow \Gamma}}$, and it follows that
\[ \bigcap_{i \in I}((\bigcap_{k \in I} ker \pi_{k_{\Omega \rightarrow
\Gamma}}) +
ker \pi_{i_{\Omega}}) \subset \bigcap_{i \in I} ker \pi_{i_{\Omega
\rightarrow
\Gamma}}. \]
Thus, $J(B)$ satisfies (\ref{adapti}). \hfill$\square$
\begin{pr} \label{diff1} 
Let $( \Gamma(B),(\Gamma(B_i))_{i \in I}) $ be
adapted to $(B,(J_i)_{i \in I})$.
Then the covering completion of $( \Gamma(B), (ker \pi_{i_\Gamma})_{i \in
I})$
is an LC differential algebra over $B_c$. 
\end{pr}
Proof: Let
\[  \pi^i_{j_{\Gamma}} : \Gamma(B)/ker\pi_{i_{\Gamma}} \longrightarrow
\Gamma(B)/ (ker \pi_{i_{\Gamma}} + ker
\pi_{j_{\Gamma}})\]
be the quotient maps.
Since the ideals $ker \pi_{i_{\Gamma}}$ are differential ideals, the
factor algebras \\
$ \Gamma(B) /ker \pi_{i_{\Gamma}}$ and
$\Gamma(B)/(ker\pi_{i_{\Gamma}}+ker\pi_{j_{\Gamma}})$ are differential
calculi
over $B_i$ and $B_{ij}$, and the projections $\pi^i_{j_{\Gamma}}$ are
the extensions of the projections $\pi^i_j$.

According to Definition \ref{C} the covering completion
$(\Gamma_c(B_c),(ker
\pi_{i_{\Gamma_c}})_{i \in I})$ of the pair\\
 $(\Gamma(B),(ker
\pi_{i_{\Gamma}})_{i \in I})$
has the entries 
\begin{eqnarray} \label{complet} 
\Gamma_c(B_c):&=&
\{(\gamma_i)_{i\in I} \in \bigoplus_i
\Gamma(B)/ ker \pi_{i_{\Gamma}}|\pi^i_{j_{\Gamma}}(\gamma_i)=
\pi^j_{i_{\Gamma}}(\gamma_j)\} \\
ker\pi_{k_{\Gamma_c}}:&=&\{ (\gamma_i)_{i\in I} \in \Gamma_c(B_c)| \gamma_k=0\}.
\end{eqnarray}

$\Gamma_c(B_c)$ has a natural grading coming from the grading of the
differential
calculi
$\Gamma(B)/ker \pi_{i _{\Gamma}}$ and a natural differential $d_c$ given by
\begin{equation} \label{dc} 
d_c(( \gamma_i)_{i\in I})=( d\gamma_i)_{i\in I} \;\; \forall
(\gamma_i)_{i\in I} \in
\Gamma_c(B_c). 
\end{equation}

From $\pi_{k_\Gamma}=\pi_{k_{\Gamma_c}}\ci K_\Gamma$ (where $K_\Gamma$
is the extension of $K:B\lra B_c$) and
formula
(\ref{dc}) it follows that
$\pi_{k_{\Gamma_c}}$ is surjective and differentiable, thus $\Gamma_c(B_c)/ker
\pi_{i_{\Gamma_c}}$ is a differential calculus isomorphic to 
$\Gamma(B)/ker\pi_{i_{\Gamma}}$.
Because of $(\Gamma(B)/ker\pi_{i_{\Gamma}})^0=B_i$, $\Gamma^0_c(B_c)=B_c$.
\hfill$\square$
\\

If the covering $(ker\pi_{i_{\Gamma}})_{i \in I}$ is complete,
$\Gamma_c(B_c)$
is isomorphic to $\Gamma(B)$ as differential calculus.
\\

Since the differential calculi
$\Gamma(B_i)$ and $\Gamma(B)/(ker \pi_{i_{\Gamma}})$ are canonically
isomorphic,
formula (\ref{complet}) is the same as
\[\Gamma_c(B) =\{(\gamma_i)_{i \in I}\in \bigoplus_i \Gamma(B_i)|
\pi^i_{j_{\Gamma}}(\gamma_i)=\pi^j_{i_{\Gamma}}(\gamma_i)\}. \]

We also note that the differential ideal $J(B_{ij}) \subset \Omega(B_{ij})$
corresponding to \\
$\Gamma(B_{ij}):= \Gamma(B)/(ker \pi_{i_{\Gamma}}+ ker \pi_{j_{\Gamma}})$ is
\begin{eqnarray*} 
J(B_{ij})&=&\pi_{{ij}_{\Omega}}(ker \pi_{i_{\Omega
\rightarrow \Gamma}} + ker \pi_{j_{\Omega \rightarrow
\Gamma}}) \\
&=&\pi^i_{j_{\Omega}}\circ \pi_{i_{\Omega}}(ker
\pi_{i_{\Omega \rightarrow
\Gamma}})+ \pi^j_{i_{\Omega}} \circ \pi_{j_{\Omega}}(ker \pi_{j_{\Omega
\rightarrow \Gamma}}), 
\end{eqnarray*} 
and because of $\pi_{i_{\Omega}}(ker
\pi_{i_{\Omega \rightarrow \Gamma}})= J(B_i)$ one can also write
\begin{equation} \label{ij} 
J(B_{ij})= \pi^i_{j_{\Omega}}(J(B_i)) +
\pi^j_{i_{\Omega}}(J(B_j)). 
\end{equation}
\\

All the above considerations remain unchanged if one considers
$*$-algebras, $*$-ideals, $*$-homomorphisms and differentials commuting ( or
anticommuting) with $*$ ($*$-differential algebras).

\section{Example}
\indent

Here we present an example of a quantum space being glued together from two
copies of a quantum disc. The result is a  $C^*$-algebra isomorphic to the
algebra of the Podle\'s spheres $S^2_{\mu c}$, $c>0$. An analogous 
construction was performed
in \cite{buko}, using another kind of quantum disc, and it was mentioned
there
in a footnote, that the resulting $C^*$-algebra is isomorphic to a Podle\'s
sphere.
To prove this isomorphy, we start with a result of Sheu \cite{sheu}, saying
that the Podle\'s spheres $S^2_{\mu c}$, $c>0$ are isomorphic as
$C^*$-algebras
to the fibered product $C^*(\mathfrak{S})\oplus_\sigma C^*(\mathfrak{S})$ of two shift
algebras by
means of the symbol map $\sigma$. In our terminology, two copies of
$C^*(\mathfrak{S})$
are glued together using the homomorphism $\sigma:C^*(\mathfrak{S})\lra C(S^1)$.
On the
other hand, using results and arguments from \cite{klile} and \cite{klile2},
it is easy to show that quantum discs are as $C^*$-algebras isomorphic to shift
algebras, and that the symbol map is transported into a natural homomorphism
of the
quantum disc onto $C(S^1)$, which just corresponds to the classical circle
contained in the quantum disc. This gives the desired isomorphy.
Moreover, using the generators of the quantum disc, we get a description  of
the glued $C^*$-algebra in terms of generators and relations. We argue that
these generators should be considered as natural ``coordinates'' on a quantum version of 
a top of a cone.
Finally, we apply our prescription of gluing together differential calculi.
Starting from $U_{q^{1/2}}(sl_2)$-covariant differential calculi on the
discs, we obtain
a differential calculus on our glued quantum top, which can also
be
characterized in terms of generators and relations, and is also
$U_{q^{1/2}}(sl_2)$-covariant.
\bde
The $C^*$-algebra $C(D_q)$, $0<q < 1$, of the quantum disc $D_q$ is defined
as
the $C^*$-closure
of the algebra $P(D_q):=\mathbbm{C}<x,x^*>/J_q$, where $J_q$ is the two-sided
ideal
in the free algebra $\mathbbm{C}<x,x^*>$ generated by the relation
\beq
x^*x-qxx^*=(1-q)1.\label{disc}
\eeq
\ede

This is a one-parameter subfamily of the two-parameter family of quantum
discs
described
in \cite{klile2}. The C*-closure is formed using only bounded *-representations
of $P(D_q)$. This is possible, because $\|\rho(x)\|=1$ for every bounded *-representation, as is shown in \cite{klile2}. From there, we also have
\bpr \label{repdq}
Every irreducible *-representation of $C(D_q)$ is unitarily equivalent to
one of the following representations:\\
(i) a one-dimensional representation $\rho_{\theta}$, defined by
$\rho_{\theta}(x)=e^{i\theta}$, $\rho_{\theta}(x^*)=e^{-i\theta}$,
for $0\leq\theta<2\pi$.\\
(ii) an infinite dimensional representation $\pi_q$ defined on a Hilbert space
$H$ with orthonormal basis $(e_i)_{i\geq 0}$ by
\beq 
\pi_q(x)e_i=\sqrt{\lambda_{i+1}}~e_{i+1},~~~i\geq 0,\label{pix}
\eeq
\beq
\pi_q(x^*)e_i=\left\{\begin{array}{lr}0&i=0,\\\sqrt{\lambda_i}~e_{i-1},&
i\geq 1,
\end{array}\right.\label{pix*}
\eeq
with $\lambda_i=1-q^i$, $i\geq 0$.
\epr

It is also shown in \cite{klile2} that the infinite dimensional representation $\pi_q$
is faithful. Therefore, $C(D_q)$ has no nontrivial covering. The
one-dimensional representations $\rho_{\theta}$ correspond to the
classical
points, forming a circle, of the quantum disc. Considering $C(S^1)$
as the $C^*$-algebra generated by $a, a^*$ with relations $aa^*=a^*a=1$, the
embedding of this classical circle into the quantum disc is described by a
$C^*$-homomorphism $\phi_q:C(D_q)\lra C(S^1)$ defined by $x\mapsto a$.
Later we will need
\ble \label{basis}
The elements $x^k{x^*}^l,~~k,l\geq 0$ form a vector space basis of $P(D_q)$.
The same is true for the elements $(xx^*)^kx^l, ~k\geq 0, ~l\in \Z$,
where $x^{-l}:={x^*}^l$, $l>0$.
\ele
Proof: It is obvious from the relations that every element of $P(D_q)$ can
be
written as a linear combination of the given elements. Applying the 
representation $\pi_q$ to an equation $\sum c_{kl}x^k{x^*}^l=0$ and acting
with
the zero operator $\pi_q(\sum c_{kl}x^k{x^*}^l)$ onto suitable basis elements
one obtains $c_{kl}=0,~\forall k,l$. The linear independence of the
second set
of elements follows from
\[(xx^*)^kx^l=q^{\frac{1}{2}k(k-1+2l)}x^{k+l}{x^*}^k+a^{kl}_{k-1}
x^{k+l-1}{x^*}^{k-1}+\ldots+a^{kl}_0x^l,~~l\geq 0,\]
\[(xx^*)^k{x^*}^l=q^{\frac{1}{2}k(k-1)}x^k{x^*}^{k+l}+
a^{k0}_{k-1}x^{k-1}{x^*}^{k-1+l}+\ldots a^{k0}_0{x^*}^l,~~l>0,\]
(with some coefficients $a^{jk}_l$), using a ``triangular type''
argument.\epf

The argument proving that the $x^k{x^*}^l$ form a basis also shows that 
$\pi_q$ is faithful on $P(D_q)$ and that consequently $P(D_q)$ is faithfully 
embedded in $C(D_q)$.
\begin{lem} \label{divisor} 
$1-xx^*$ is not a zero divisor in $P(D_q)$.
\end{lem}
Proof: Assume $(1-xx^*)\sum_{kl} c_{kl} (xx^*)^k x^l=0$. Lemma \ref{basis}
 gives the
following conditions for the coefficients $c_{kl}$,
\begin{eqnarray*} 
c_{0l}&=&0,~~ \forall l \\
c_{kl}&=&c_{k-1,l},~~\forall l,k \geq 1 
\end{eqnarray*} 
which lead to
$c_{kl}=0,~~ \forall k,l$. In
the same way one proves that $1-xx^*$ is also not a right zero divisor.
\hfill$\square$
\\

The $C^*$-algebra of the unilateral shift is defined as follows:
Let $H$ be a Hilbert space with orthonormal basis $(e_i)_{i\geq 0}$. The
shift operator ${\mathfrak{S}} \in B(H)$ is defined by $\mathfrak{S}(e_i)=e_{i+1}$.
Its adjoint is given
by
$\mathfrak{S}^*(e_i)=\left\{\begin{array}{lr}e_{i-1}&i>0\\0&i=0.\end{array}\right.
$
$C^*(\mathfrak{S})$ is the $C^*$-subalgebra of $B(H)$ generated by 
$\mathfrak{S}$ and
$\mathfrak{S}^*$. The symbol
map $\sigma:C^*(\mathfrak{S})\lra C(S^1)$ is the homomorphism defined by
$\sigma(\mathfrak{S})=a$.
\bpr\label{dishi}
$C(D_q)$ is isomorphic to $C^*(\mathfrak{S})$ as $C^*$-algebra. Under this
isomorphism,
the symbol map $\sigma$ corresponds to the embedding of the classical
circle,
$\phi_q:C(D_q)\lra C(S^1)$.
\epr
Proof: We use ideas of \cite{klile}, where this is proved for another
one-parameter
family of quantum discs. In fact, we prove $\pi_q(C(D_q))=C^*(\mathfrak{S})$. First,
it is easy to see that
\beq
\pi_q(x)=\mathfrak{S}\sum_{k=0}^\infty (\sqrt{\lambda_{k+1}}-\sqrt{\lambda_k})
\mathfrak{S}^k{\mathfrak{S}^*}^k,\label{piS}
\eeq
where the series converges in the operator norm. Thus, $\pi_q(C(D_q))\subset
C^*(\mathfrak{S})$. On the other hand, if $P_i$ is the orthogonal projector onto
$e_i$,
$\pi_q(x^*)\pi_q(x)=\sum \lambda_{i+1}P_i$ is the spectral resolution of
$\pi_q(x^*)\pi_q(x)$, and $P_i$ lies in the $C^*$-algebra
generated by $\pi_q(x^*)\pi_q(x)$, therefore also in $\pi_q(C(D_q))$. The matrix
units $E_{ij}$ defined by $E_{ij}(e_k)=\delta_{jk}e_i$ can be written
\[E_{ij}=\left\{\begin{array}{lr}
(\sqrt{\lambda_{j+1}\cdots\lambda_i})^{-1} \pi_q(x)^{i-j}P_j &i>j,\\
(\sqrt{\lambda_j\cdots\lambda_{i+1}})^{-1} \pi_q(x^*)^{j-i}P_j & i<j, \\
P_i&i=j. \end{array}\right.\]

Since the $E_{ij}$ generate the ideal ${\cal K}$ of compact operators, it
follows that ${\cal K}\subset\pi_q(C(D_q))$. Moreover, $\mathfrak{S}-\pi_q(x)$ is a
weighted
shift, $(\mathfrak{S}-\pi_q(x))(e_k)=(1-\sqrt{\lambda_{k+1}})e_{k+1}$, where
$1-\sqrt{\lambda_k}\lra 0$
for $k\lra \infty$. From the next lemma it follows that $\mathfrak{S}-\pi_q(x)\in
{\cal K}$,
therefore $\mathfrak{S} \in\pi_q(C(D_q))$.
For $\sigma\ci\pi_q=\phi_q$ it is sufficient to show
$\sigma(\pi_q(x))=\phi_q(x)=a$,
which follows from formula (\ref{piS}) using
$\sum_{i=0}^\infty(\sqrt{\lambda_{i+1}}-
\sqrt{\lambda_i})=\lim_{k\to\infty}\sqrt{\lambda_k}=1$.
\epf
\ble
Let $T\in B(H)$ be a weighted shift,
\[T(e_i)=t_ie_{i+1},\]
with $t_i\in \mathbbm{R}$, $\lim_{i\to\infty}t_i=0$. Then $T\in {\cal K}$.
\ele
Proof: $T^*(e_i)=\left\{\begin{array}{lr}0&i=0,\\t_{i-1}e_{i-1}&i>0.
\end{array}\right.$ Therefore, $T^*T(e_i)=t_i^2e_i$, and $T^*T$ is a compact
operator with spectrum consisting of the isolated eigenvalues $t_i^2$,
$t_i^2\to 0$. Then $\sqrt{T^*T}$ is also a compact operator with eigenvalues
$|t_i|$, and $T$ itself is compact, because its polar decomposition is
$T=U\sqrt{T^*T}$, and ${\cal K}$ is a two-sided ideal.\epf
\\

Now, the following proposition is immediate from Proposition 1.2. of
\cite{sheu}.
\bpr
For $|\mu|<1,~c>0$, the $C^*$-algebra $C(S^2_{\mu c})$ of the Podle\'s sphere
is
isomorphic to $C(D_p)\oplus_\phi C(D_q)=\{(f,g)\in C(D_p)\oplus C(D_q)|
\phi_p(f)=\phi_q(g)\}$, $0<q,p < 1$.
\epr

The isomorphy holds for any pairs of parameters $(\mu,c)$ and $(p,q)$.
The images of the generators of $C(S^2_{\mu c})$ under this isomorphism clearly
are also generators of $C(D_p)\oplus_\phi C(D_q)$. However, one can also
describe $C(D_p)\oplus_\phi C(D_q)$ by means of generators arising naturally
from the generators of the two quantum discs via the gluing procedure:
%
\bpr \label{pqtop}
Let
\[P(S^2_{pq\phi}):=\mathbbm{C}<f_1,f_{-1},f_0>/J_{q,p},\]
where $J_{q,p}$ is the two-sided ideal generated by the relations
\beq
f_{-1}f_1-qf_1f_{-1}=(p-q)f_0+(1-p)1,\label{re1}
\eeq
\beq
f_0f_1-pf_1f_0=(1-p)f_1,\label{re2}
\eeq
\beq
f_{-1}f_0-pf_0f_{-1}=(1-p)f_{-1},\label{re3}
\eeq
\beq
(1-f_0)(f_1f_{-1}-f_0)=0.\label {re4}
\eeq
With 
\begin{equation} \label{re*} 
f_0^*=f_0,~~f_1^*=f_{-1}, 
\end{equation}
$P(S^2_{pq\phi})$ is a $*$-algebra.\\
$P(D_p)\oplus_\phi P(D_q)=\{(f,g)\in P(D_p)\oplus P(D_q)|
\phi_p(f)=\phi_q(g)\}$ is isomorphic to 
the *-algebra $P(S^2_{pq\phi})$.
\epr
Remark: We have the conjecture that $P(S^2_{pq\phi})$
and $P(S^2_{\mu c})$ are not isomorphic as *-algebras.\\[.3cm]
Proof: Relation (\ref{re1}) is invariant under $*$ whereas (\ref{re2}) and
(\ref{re3}) are transformed into each other. (\ref{re4}) is invariant
because of
$f_1 f_{-1}f_0=f_0f_1f_{-1}$ which follows from (\ref{re2}) and (\ref{re3}).
\\Next we show $P(D_p)\oplus_\phi P(D_q)\simeq P(S^2_{pq\phi})$. We denote
the
generators of $P(D_p)$ by $x,x^*$ and those of $P(D_q)$ by $y,y^*$.
Consider the elements
$\tilde{f}_0=(xx^*,1),~\tilde{f}_1=(x,y),~\tilde{f}_{-1}
=(x^*,y^*)$ of $P(D_p)\oplus_\phi P(D_q)$. They fulfill the relations
(\ref{re1}) - (\ref{re4}).
\ble
$\tilde{f}_0,~\tilde{f}_1,~\tilde{f}_{-1}$ generate $P(D_p)\oplus_\phi
P(D_q)$.
\ele
Proof of the lemma: Use the basis $(xx^*)^kx^l,~(yy^*)^k y^l,~~k\geq 0,~l\in
\Z$ of
Lemma \ref{basis}. First we notice that
\[(\sum_{m\geq 0,n\in \Z}c_{mn}(xx^*)^mx^n,\sum_{k\geq 0,l\in \Z}
\tilde{c}_{kl}(yy^*)^ky^l)\in P(D_p)\oplus_\phi P(D_q) \] if and only if
\beq
\sum_{k\geq 0}c_{kl}=\sum_{k\geq 0}\tilde{c}_{kl},~~\forall l.\label{cglu}
\eeq
We obtain
\begin{eqnarray*} 
&& (\sum_{m\geq 0,n\in \Z}c_{mn}(xx^*)^mx^n,\sum_{k\geq
0,l\in \Z}
\tilde{c}_{kl}(yy^*)^ky^l)=\\ &&\sum_{l\in \Z}(\sum_{m\geq 0}c_{ml}
(xx^*)^mx^l,\sum_{k\geq 0}
\tilde{c}_{kl}(yy^*)^ky^l)=\\ &&
\sum_{l\in \Z}((\sum_{m\geq 0}c_{ml}(xx^*)^mx^l,\sum_{k\geq 0}
\tilde{c}_{kl}y^l)+\\ &&
(\sum_{m\geq 0}c_{ml}x^l,\sum_{k\geq 0}
\tilde{c}_{kl}(yy^*)^ky^l)-(\sum_{m\geq 0}c_{ml}x^l,\sum_{k\geq 0}
\tilde{c}_{kl}y^l))=\\ &&
\sum_{l\in \Z}(\sum_{m\geq 0}c_{ml}{\tilde{f}_0}^m{\tilde{f}_1}^l
+\sum_{k\geq
0}\tilde{c}_{kl}(\tilde{f}_1\tilde{f}_{-1}-\tilde{f}_0+1)^k{\tilde{f}_1}^l
-\sum_{n\geq 0}c_{nl}{\tilde{f}_1}^l), 
\end{eqnarray*}
where we have set ${\tilde{f}_1}^{-1}=\tilde{f}_{-1}$, and (\ref{cglu}) has
been used in the last equality. {}\hfill$\square$
\\

By the lemma there exists a surjective homomorphism $F: P(S^2_{pq\phi})
\longrightarrow P(D_p) \bigoplus_{\phi} P(D_q)$ defined by
$F(f_i):=\tilde{f}_i$.
Let $ p_1: P(D_p) \bigoplus_{\phi} P(D_q) \longrightarrow P(D_p)$ and
$p_2:
P(D_p) \bigoplus_{\phi} P(D_q) \longrightarrow P(D_q)$ be the restrictions
of
the first and second projections.
By definition $ker\,p_1 \bigcap ker\, p_2=0$. Let $\pi_1:=p_1 \circ F$
and $\pi_2:= p_2 \circ F$, i. e. $\pi_1(f_0)=xx^*, ~\pi_1(f_1)=x,
\,\,\pi_2(f_0)=1,~\pi_2(f_1)=y$. $F$ is an isomorphism if $ker \pi_1 \bigcap
ker \pi_2 =0$.

First we describe $ker \pi_2$. Every element $a \in P(S^2_{pq\phi})$ can be
written in the form \[ a= \sum_{k\in {\Z}; m,n \geq 0} c_{mnk} (f_1
f_{-1})^m f_0^n
f_1^k, \] where $f^{-1}_1=f_{-1}$. Applying $\pi_2$ to $a \in ker\,\pi_2$,
it
follows that
\[ \sum_{k \in {\Z};m,n \geq 0} c_{mnk} (yy^*)^m y^k =0, \] and one obtains
the
condition $\sum_{n \geq 0} c_{mnk}=0,~~\forall m,k$.
Thus, we have the identity 
\begin{eqnarray*} 
\sum_{n \geq 0} c_{mnk}
(f_1f_{-1})^m f^n_0
f^k_1 &=& \sum_{n \geq 1} \sum^n_{s= 1} c_{mnk} (f_1 f_{-1})^m (f^s_0 -
f_0^{s-1})f^k_1\\ &=&\sum_{n \geq 1} \sum^n_{s=1} c_{mnk} (f_1f_{-1})^m
(f_0-1)f^{s-1}_0 f^k_1. 
\end{eqnarray*}
In view of $f_1f_{-1}f_0=f_0f_1f_{-1}$, this means that every element $a \in
ker \pi_2$ can be written in the form
\begin{equation}  \label{acoeff} 
a=(1 -f_0)\sum_{k \in {\Z}; m,n \geq 0}
a_{mnk} (f_1f_{-1})^m
f^n_0 f^k_1. 
\end{equation}

Assume now $a \in ker \pi_1 \bigcap ker \pi_2$. Applying $\pi_1$ to $a$
one obtains 
\[ (1-xx^*)\sum_{k \in {\Z}; m,n  \geq 0} a_{mnk} (xx^*)^m (xx^*)^n x^k=0, \] which
yields, since
$1-xx^*$ is not a zero divisor, the following condition for the
coefficients $a_{mnk}$: 
\begin{equation} \label{1minus} 
\sum^l_{n = 0}
a_{l-n,n,k}=0,~\forall l \geq 0,k. 
\end{equation}
This leads to 
\begin{eqnarray*} && (1-f_0)\sum_{k \in {\Z}, m,n \geq 0}
a_{mnk}(f_1
f_{-1})^m f^n_0f^k_n \\&=& (1-f_0) \sum_{k \in {\Z}} \sum_{l \geq 0}
\sum^l_{n=0} a_{l-n,n,k} (f_1f_{-1})^{l-n}
f^n_0 f^k_n\\ &=& (1-f_0) \sum_{k \in {\Z}} \sum_{l \geq 0} \sum^l_{n=1}
\sum^n_{s=1} a_{l-n,n,k}
((f_1f_{-1})^{l-s}f^s_0-(f_1f_{-1})^{l-s+1} f_0^{s-1})f^k_1 
\end{eqnarray*}
\begin{equation} \label{kerpi1} 
=-(1-f_0)(f_1f_{-1}-f_0)\sum_{k \in {\Z}}
\sum_{l \geq 0} \sum^l_{n=1}\sum^n_{s=1} a_{l-n,n,k}
(f_1f_{-1})^{l-s} f^{s-1}_0 f^k_1=0. 
\end{equation}
(\ref{1minus}) was used in the second equality.
Thus $ker\, \pi_1 \bigcap ker \,\pi_2=0$, and $F$ is an isomorphism.\epf
\vspace{.3cm}

Note that the computation leading to (\ref{kerpi1}) can also be used to
show
that $ker \pi_1$ is generated by $f_1f_{-1}-f_0$, since the application of
$\pi_1$ to a general element of $ker \pi_1$ leads to (\ref{kerpi1}) without
the
factor $1-f_0$.
\\\\
\bpr
Let $\rho$ be a homomorphism of the *-algebra $P(S^2_{pq\phi})$ into the
*-algebra $B(H)$ of bounded operators on a Hilbert space $H$.\\
Then $ker~(1-\rho(f_0))$ and $ker~\rho(f_1f_{-1}-f_0)$ are closed 
subspaces invariant under all representation operators.
$H$ can be decomposed into the orthogonal direct sum
\beq\label{kerzer}
H=ker(1-\rho(f_0))\oplus ker(\rho(f_1f_{-1}-f_0)|_{(ker(1-\rho(f_0)))^\bot}).
\eeq
\epr
Proof: The invariance of the kernels is a direct consequence of the relations.
The kernels are closed since they belong to bounded operators. Thus, $H$ can
be decomposed into the orthogonal direct sum $H=ker(1-\rho(f_0))\oplus 
ker(1-\rho(f_0))^\bot$ of closed invariant subspaces. In turn,
$ker(1-\rho(f_0))^\bot$ can be decomposed as
$ker(\rho(f_1f_{-1}-f_0)|_{(ker(1-\rho(f_0)))^\bot})\oplus H^c$, with another
invariant subspace $H^c$. The restrictions to $H^c$ of both $1-\rho(f_0)$
and $\rho(f_1f_{-1}-f_0)$ are injective. Thus, if $\psi\in H^c, \psi\neq 0$,
it follows that $(1-\rho(f_0))\rho(f_1f_{-1}-f_0)(\psi)\neq 0$, which 
contradicts relation (\ref{re4}).\epf

In the restriction of $\rho$ to $ker(1-\rho(f_0))$ we have $\rho(f_0)=1$, 
and the relations reduce to 
\[\rho(f_{-1})\rho(f_1)-q\rho(f_1)\rho(f_{-1})=(1-q)1,\]
which are the relations of a quantum disc with parameter $q$. On the complement
$ker(\rho(f_1f_{-1}-f_0)|_{(ker(1-\rho(f_0)))^\bot})$ we have 
$\rho(f_0)=\rho(f_1f_{-1})$, and the relation (\ref{re1}) reduces to
\[\rho(f_{-1})\rho(f_1)-p\rho(f_1)\rho(f_{-1})=(1-p)1,\]
i. e. the relations of a quantum disc with parameter $p$. (\ref{re2})
and (\ref{re3}) follow from these relations, whereas (\ref{re4}) is satisfied
trivially. Using the results of \cite{klile2} for quantum discs, we obtain
\bpr
The following
is a complete list of bounded irreducible *-representations of
$P(S^2_{pq\phi})$:\\
1. A representation in a Hilbert space $H$ with orthonormal basis
$(e_i)_{i=0,1,...}$
\beq
\rho_1(f_0)e_i=\left\{\begin{array}{lr}0&i=0\\ \lambda_i e_i&i>0\end{array}
\right.
\label{pi1f0}
\eeq
\beq
\rho_1(f_1)e_i=\sqrt{\lambda_{i+1}}~e_{i+1},~~~i\geq 0,\label{pi1f1}
\eeq
\beq
\rho_1(f_{-1})e_i=\left\{\begin{array}{lr}0&i=0,\\\sqrt{\lambda_i}~e_{i-1},&
i\geq 1,
\end{array}\right.\label{pi1f-1}
\eeq
with $\lambda_i=1-p^i$, $i\geq 0$. 
\\
2. A representation in $H$ given by
\beq
\rho_2(f_0)e_i=e_i,\label{pi2f0}
\eeq
\beq
\rho_2(f_1)e_i=\sqrt{\lambda'_{i+1}}~e_{i+1},~~~i\geq 0,\label{pi2f1}
\eeq
\beq
\rho_2(f_{-1})e_i=\left\{\begin{array}{lr}0&i=0,\\\sqrt{\lambda'_i}~e_{i-1},
&i\geq 1,
\end{array}\right.\label{pi2f-1}
\eeq
with $\lambda'_i=1-q^i$, $i\geq 0$.
\\
3. A one parameter family of one dimensional representations given by
\beq
\rho_\theta(f_0)=1,\label{rof0}
\eeq
\beq
\rho_\theta(f_1)=e^{i\theta},\label{rof1}
\eeq
\beq
\rho_\theta(f_{-1})=e^{-i\theta},\label{rof-1}
\eeq
where $0\leq\theta<2\pi$.\\
Moreover, $\|\rho(f_1)\|=\|\rho(f_{-1})\|=\|\rho(f_0)\|=1$ for any
*-representation of $P(S^2_{pq\phi})$ in bounded operators.
\epr

Denoting by $Rep_b$ the set of *-representations of $P(S^2_{pq\phi})$ in
bounded operators, it follows that for each $a\in P(S^2_{pq\phi})$ exists
$\|a\|:=\sup_{\rho\in Rep_b}\|\rho(a)\|<\infty$. 
\bpr\label{basps2}
(i) $\|.\|$ is a C*-norm on $P(S^2_{pq\phi})$.\\
(ii) $\rho_1\oplus\rho_2$ is a faithful representation of $P(S^2_{pq\phi})$.\\
(iii) $\{f_1^kf_0f_{-1}^l, f_1^kf_{-1}^l|k,l=0,1,\ldots\}$ is a vector space
basis of $P(S^2_{pq\phi})$.
\epr
Proof: As a first step, one shows that the vectors (iii) form a linear
generating system. For this, one first shows inductively that the monomials
$f_1^kf_0^lf_{-1}^m,~k,l,m=0,1,\ldots$ form a linear generating system.
Then one uses (\ref{re4}) to reduce the power of $f_0$.\\
Let now $a=\sum_{k,l=0,1\ldots}(a_{kl}f_1^kf_0f_{-1}^l+b_{kl}f_1^kf_{-1}^l)$,
and assume $\rho_1\oplus\rho_2(a)=0$. Then we have
\[\rho_1(a)e_0=\sum_k(a_{k0}\rho_1(f_1)^k\rho_1(f_0)e_0+b_{k0}\rho_1(f_1)^ke_0)
=\sum_kb_{k0}\sqrt{\lambda_k\cdots\lambda_1}e_k=0,\]
i. e. $b_{k0}=0$ for all $k$. Thus,
\[\rho_2(a)e_0=\sum_k(a_{k0}\sqrt{\lambda_k'\cdots\lambda_1'}e_k+
b_{k0}\sqrt{\lambda_k'
\cdots\lambda_1'}e_k)=0\]
gives $a_{k0}=0$ for all $k$.

Assume now that
\[a_{kl}=b_{kl}=0,~ \forall k\]
is shown for $l\leq i$. Then
\[\rho_1(a)e_{i+1}=\sum_k(a_{k,i+1}\rho_1(f_1)^k\rho_1(f_0)\sqrt{\lambda_1\cdots
\lambda_{i+1}}e_0+b_{k,i+1}\rho_1(f_1)^k\sqrt{\lambda_1\cdots\lambda_{i+1}}e_0)
=\]
\[=\sum_kb_{k,i+1}\sqrt{\lambda_k\cdots\lambda_1^2\cdots \lambda_{i+1}}e_k=0,\]
i. e. $b_{k,i+1}=0$. $\forall k$. From 
\[\rho_2(a)e_{i+1}=\sum_k(a_{k,i+1}\sqrt{\lambda_k'\cdots{\lambda_1'}^2\cdots
\lambda_{i+1}'}e_k+b_{k,i+1}\sqrt{\lambda_k'
\cdots{\lambda_1'}^2\cdots\lambda_{i+1}'}e_k)=0\]
now also follows $a_{k,i+1}=0$, $\forall k$. This proves the proposition.\epf
\bde
$C(S^2_{pq\phi})$ is the closure of $P(S^2_{pq\phi})$ in the norm $\|.\|$.
\ede
%
\bpr
$C(S^2_{pq\phi})$ is C*-isomorphic to $C(S^2_{\mu c})$ for $c>0$.
\epr
Proof: It is sufficient to prove that $C(S^2_{pq\phi})$ is isomorphic to
$C^*(\mathfrak{S})\oplus_\sigma C^*(\mathfrak{S})$. Using 
that reduced atomic representations are faithful, it is enough to see
that one can choose for both algebras one element from every equivalence class
of irreducible representations in such a way that there is a bijection
between the resulting sets of representations, and that the images of 
corresponding representations are isomorphic as C*-algebras.
It follows again from Proposition
1.2 of
\cite{sheu} that the irreducible representations of
$C^*(\mathfrak{S})\oplus_\sigma C^*(\mathfrak{S})$
are up to equivalence $p'_1$, $p'_2$ and $p_\theta\ci\sigma$, where 
$p'_1$, $p'_2$ are the
restrictions of the first and second projections to 
$C^*(\mathfrak{S})\oplus_\sigma C^*(\mathfrak{S})$, and $p_{\theta}$ is the
evaluation at
$e^{i\theta}$, $p_\theta(f)=f(e^{i\theta})$ for $f\in C(S^1)$.
Indeed, $p'_1$, $p'_2$, $p_\theta \circ \sigma$ correspond under the
isomorphism
$(\pi_+,\pi_-):C(S^2_{\mu c})\lra C^*(\mathfrak{S})\oplus_\sigma C^*(\mathfrak{S})$
to the
representations $\pi_+,\pi_-,\pi_\theta$ of $C(S^2_{\mu c})$ (see \cite{po})
respectively. The representations of $C(S^2_{pq\phi})$ corresponding to
$p'_1$, $p'_2$, $p_\theta\ci\sigma$ are now $\rho_1$, $\rho_2$, $\rho_\theta$
respectively. The equality of the corresponding images is trivial for the
one-dimensional representations and follows for the others with the same 
arguments as in the proof of Proposition \ref{dishi}.\epf

Note that it would have been sufficient to use only $\rho_1, \rho_2$ and
$p'_1,p'_2$ in the above proof, because $\rho_1\oplus \rho_2$ and $p'_1
\oplus p'_2=id$ are already faithful representations. For $\rho_1\oplus
\rho_2$ this follows from $\rho_\theta=p_\theta\ci\sigma\ci\rho_i,~i=1,2$.
Moreover, there are the equalities
$p'_i\ci(\pi_p\oplus\pi_q)\ci F=\rho_i$, $i=1.2$, and $p_\theta\ci\sigma
\ci(\pi_p\oplus\pi_q)\ci F=\rho_\theta$, which mean that $F$ extends to
a C*-isomorphism $C(S^2_{pq\phi})\lra C(D_p)\oplus_\phi C(D_q)$.
\\

In order to determine an underlying ``space''
of $C(S^2_{pq\phi})$, we look for the spectra of generators. First we
introduce instead of $f_1,f_{-1}$ the selfadjoint elements
$f_+=\frac{1}{2}(f_1+f_{-1})$, $f_-=\frac{1}{2}i(f_1-f_{-1})$.
In terms of $f_+$ and $f_-$ the relations (\ref{re1}) - (\ref{re4}) are
\beq
(1-q)(f_+^2 + f_-^2) + (1+q)i(f_-f_+-f_+f_-)=(p-q)f_0+(1-p)1,
\label{re1+-}\eeq
\beq
f_0f_+ -pf_+f_0-i(f_0f_--pf_-f_0)=(1-p)(f_+-if_-),
\label{re2+-}\eeq
\beq
f_+f_0-pf_0f_++i(f_-f_0-pf_0f_-)=(1-p)(f_++if_-),\label{re3+-}
\eeq
\beq
(1-f_0)(f_+^2+f_-^2+i(f_+f_--f_-f_+)-f_0)=0.
\label{re4+-}\eeq
Putting here $p=q=1$, (\ref{re1+-}) - (\ref{re3+-}) just mean
commutativity of
the algebra, whereas the geometric counterpart of (\ref{re4+-}) is
the union of the plane $f_0=1$ and the cone $f_+^2+f_-^2=f_0$ in
$f_0,f_+,f_-$-space.
For $p,q\neq 1$, $f_+$ and $f_-$ act in the irreducible representations as
follows :
\beq
\rho_1(f_+)e_k=\frac{1}{2}(\sqrt{\lambda_{k+1}}e_{k+1}+\sqrt{\lambda_k}e_{k-
1}),\label{r1f+}
\eeq
\beq
\rho_1(f_-)e_k=\frac{1}{2}i(\sqrt{\lambda_{k+1}}e_{k+1}-\sqrt{\lambda_k}e_{k
-1}),\label{r1f-}
\eeq
$\rho_2(f_\pm)$ obey (\ref{r1f+}) and (\ref{r1f-}) with $\lambda_k'$ in
place of
$\lambda_k$, and
\beq
\rho_\theta(f_+)=\cos\theta,~~~\rho_\theta(f_-)=\sin\theta.
\eeq
It follows that
\[\rho_1(f_+^2+f_-^2)e_i=(1-\frac{1}{2}(p^i+p^{i+1}))e_i\]
and
\[\rho_2(f_+^2+f_-^2)e_i=(1-\frac{1}{2}(q^i+q^{i+1}))e_i,\]
whereas $\rho_{1,2}(f_\pm)$ are Jacobi ``matrices'' with continuous spectra. So
one
may draw the following picture in $f_0,f_+,f_-$-space, assigning to 
every ``eigenstate''
$e_i$ of  $\rho_{1,2,\theta}(f_0)$ the possible values of a
``measurement'' of $f_0, f_+,f_-$:\\
$\rho_\theta$, $0\leq\theta<2\pi$ give rise to the circle
$f_0=1$, $f_+^2+f_-^2=1$. $\rho_2$ leads to circles $f_0=1$, $f_+^2+f_-^2=
1-\frac{1}{2}(q^i+q^{i+1})$, and $\rho_1$ to circles $f_0=1-p^i$,
$f_+^2+f_-^2=1-\frac{1}{2}(p^i+p^{i+1})$. The union of all these circles
may be considered as a discretized version of the top of the cone arising
in the
classical case, i. e. the set $\{(f_0,f_+,f_-)|f_0=1, f_+^2+f_-^2 \leq
1\}\cup
\{(f_0,f_+,f_-)|0\leq f_0\leq 1, f_+^2+f_-^2=f_0\}$.
\\

As a consequence of the considerations in the proof of Proposition 
\ref{pqtop} the homomorphisms $\pi_{1,2}:P(S^2_{pq\phi})\lra P(D_{p,q})$
define a covering $(ker\,\pi_1,ker\,\pi_2)$, which is just the transport
of the covering $(ker\,p_1,ker\,p_2)$ of $P(D_p)\oplus_\phi P(D_q)$
under the isomorphism $F$.\vspace{.5cm}

Our task is now the construction of a differential calculus over
$P(S^2_{pq\phi})$ adapted to the covering $(ker\,\pi_1, ker\, \pi_2)$.
According to Proposition \ref{diff}, such a differential calculus is
uniquely
determined by differential calculi over $P(D_p) \simeq
P(S^2_{pq\phi})/ker\,
\pi_1$ and $P(D_q) \simeq P(S^2_{pq\phi})/ker\,\pi_2$. \\We choose
$\Gamma(P(D_p))=\Omega(P(D_p))/J(P(D_p))$, where $J(P(D_p))$ is generated by
the
elements:
\begin{eqnarray*} x(dx) &-& p^{-1}(dx)x \\ x^*(dx^*)&-& p(dx^*)x^* \\
x(dx^*) &-& p^{-1}(dx^*)x \\ x^*(dx) &-& p(dx)x^*, \end{eqnarray*}
analogously
for $P(D_q)$. These differential calculi were considered in \cite{vaks}.
Obviously, they are $*$-differential calculi.
\begin{lem} \label{modbas} (i) $dx$ and $dx^*$ form a left and right
$P(D_q)$-module basis
of
$\Gamma^1(P(D_q))$. \\
(ii) $dxdx^*$ is a left and right $P(D_q)$-module basis of 
$\Gamma^2(P(D_q))$.\\
(iii) $\Gamma^n(P(D_q))=0$ for $n\geq 3$.
\end{lem}
This is proved in the appendix.
\\

For $p \not= q$ the differential calculus $\Gamma(P(S^2_{pq\phi}))$
obtained by gluing together $\Gamma(P(D_p))$ and $\Gamma(P(D_q))$
has the following form:
\begin{eqnarray*} \Gamma^0(P(S^2_{pq\phi})) &=&
P(S^2_{pq\phi}) \\ \Gamma^n(P(S^2_{pq\phi}))&=&
\Gamma^n(P(D_q)) \bigoplus \Gamma^n(P(D_p)),~~\forall n>0. \end{eqnarray*}
This follows from formula (\ref{ij}), which in our case reads
$J(P(S^1))=\phi_{p_{\Omega}}(J(P(D_p)) + \phi_{q_{\Omega}}(J(P(D_q))$.
We get in the differential
ideal $J(P(S^1))$ elements of the form \begin{eqnarray*} a(da) -
p^{-1}(da)a &=& \phi_{p_{\Omega}}(x(dx) - p^{-1}(dx)x)\\
a(da) - q^{-1} (da)a &=& \phi_{q_{\Omega}}(y(dy)-q^{-1}(dy)y),
\end{eqnarray*}
where $a$ is the generator of $P(S^1)$, which leads to
$(q^{-1}-p^{-1})da
\in J(P(S^1))$ and $da=0$ in $\Gamma(P(S^1))$. In the same way follows
$da^*=0$ in $\Gamma(P(S^1))$, which means $\Gamma^n(P(S^1))=0,~ \forall
n>0$. So, there is no gluing in all degrees $n >0$.

In the case $q=p$ the differential calculus on
$P(S^2_{pq\phi})$ in higher degree than zero is not simply the
direct sum of the differential calculi on the quantum discs.
The differential ideal $J(P(S^2_{qq\phi}))$
defining the differential calculus $\Gamma(P(S^2_{qq\phi}))$
can be written in terms of the generators $f_1$, $f_{-1}$ and $f_0$ as
follows:
\begin{pr} Let $\Gamma(P(S^2_{qq\phi})):=
\Omega(P(S^2_{qq\phi}))/J(P(S^2_{qq\phi}))$,
where the
differential ideal \\
$J(P(S^2_{qq\phi}))$ is generated by the elements
\begin{eqnarray}
&& f_1(df_1) - q^{-1}(df_1)f_1,~~ f_{-1}(df_{-1}) - q(df_{-1})f_{-1},
\label{refd1}\\
&& f_1(df_{-1}) - q^{-1} (df_{-1})f_1,~~ f_{-1}(df_1) - q(df_1)f_{-1},
\label{refd2}\\
&& f_0(df_1) - (df_1)f_0,~~ f_0(df_{-1})-(df_{-1})f_0  \label{refd3}\\
&& (df_0)(f_1f_{-1} - f_0),~~(f_1f_{-1}-f_0)df_0,
\label{refd4}\end{eqnarray} and
\begin{eqnarray} \label{refd21} (1-q)df_0 df_{-1} -qf_{-1}df_0df_0 \\
\label{refd22} (1-q) df_1 df_0 -q f_1df_0df_0 \\
\label{refd23}(1-f_0)((1-q)df_1 df_{-1} -df_0 df_0) \\ \label{refd24}
(f_1f_{-1}-f_0) df_0 df_0. \end{eqnarray}
Then the homomorphisms $\pi_1$ and $\pi_2$ are differentiable and
\begin{equation} \label{ker0} ker \pi_{1_{\Gamma}} \bigcap ker
\pi_{2_{\Gamma}}=\{0\}, \end{equation}
i.e. $\Gamma(P(S^2_{qq\phi}))$ is the unique differential calculus such
that \\$(\Gamma(P(S^2_{qq\phi})), (\Gamma(P(D_q)), \Gamma(P(D_q)) ))$ is
adapted to $(P(S^2_{qq\phi}),(ker\,\pi_1, ker\,\pi_2))$ according to 
Proposition \ref{diff}. \end{pr}
Proof: One shows easily that $\pi_{1_{\Omega}}(J(P(S^2_{qq\phi}))) \subset
J(P(D_q))$ and $\pi_{2_{\Omega}}(J(P(S^2_{qq\phi}))) \subset J(P(D_q))$,
i.e. $\pi_1$ and $\pi_2$ are differentiable.

First let us prove the assertion (\ref{ker0}) for the first degree
$\Gamma^1(P(S^2_{qq\phi}))$,
i. e. $(ker \pi_{1_{\Gamma}})^1 \bigcap (ker \pi_{2_{\Gamma}})^1=\{0\}$.
Using Proposition \ref{basps2}, (iii), and (\ref{refd1}), (\ref{refd2}) and
(\ref{refd3}) one 
finds that every element $\gamma \in \Gamma(P(S^2_{qq\phi}))$ can be written in
the form 
\begin{eqnarray*} 
\gamma &=& \sum_{k,l \geq 0} (a_{k,l} f^k_1 f_0 f^l_{-1} + 
b_{k,l} f^k_1 f^l_{-1}) df_1 \\
&+& \sum_{k,l \geq 0}(c_{k,l} f^k_1 f_0 f^l_{-1}  + d_{k,l} 
f^k_1f^l_{-1})df_{-1} \\ 
&+& \sum_{k,l \geq 0}(e_{k,l} f^k_1 f_0 f^l_{-1}  + g_{k,l}
f^k_1f^l_{-1})df_0. \end{eqnarray*}
Assuming $\gamma \in ker \pi_{2_{\Gamma}}$ one obtains 
\begin{eqnarray*} \pi_{2_{\Gamma}}(\gamma)&=&  
\sum_{k,l \geq 0} (a_{k,l} 
y^k{y^*}^l +b_{k,l}y^k{y^*}^l)dy 
\\&+& \sum_{k,l \geq 0} (c_{k,l}
y^k{y^*}^l +d_{k,l}y^k{y^*}^l)dy^*=0. 
\end{eqnarray*}

Using the bimodule basis $\{dy,dy^*\}$ of $\Gamma(P(D_q))$ and the
vector space basis $\{y^k{y^*}^l|k,l=0,1,\ldots\}$ of $P(D_q)$ 
one obtains 
$a_{k,l}=-b_{k,l}$ and $c_{k,l}=- d_{k,l}$. It follows that an element 
$\gamma\in ker \pi_{2_{\Gamma}}$ can be written in the form
\begin{eqnarray*}
\gamma &=& \sum_{k,l \geq 0} a_{k,l} f^k_1 (f_0-1) f^l_{-1} df_1 \\
&+& \sum_{k,l \geq 0}c_{k,l} f^k_1 (f_0-1) f^l_{-1}df_{-1} \\
&+& \sum_{k,l \geq 0}(e_{k,l} f^k_1 f_0 f^l_{-1}  + g_{k,l}
f^k_1f^l_{-1})df_0. 
\end{eqnarray*}
The relations in the algebra give the identities 
\begin{eqnarray*}
(f_0 -1)f_i &=& q^if_i (f_0-1),~i=-1,0,1, \\
f^n_{-1}f_1&=& q^n f_1 f^n_{-1} + (1-q^n)f^{n-1}_{-1},~n>0. 
\end{eqnarray*}
From (\ref{refd4}) follows 
\begin{eqnarray*}
f_1f_{-1} df_0&=&f_0df_0, \\
(f_0-1) f_{-1} df_1 &=&q (f_0-1)df_0 - q(f_0-1)f_1 df_{-1}. 
\end{eqnarray*}
The last four equations, together with (\ref{re3}), (\ref{re4}) and 
(\ref{refd2}) now  give the following two identities:
\begin{eqnarray*} 
f_0 f^l_{-1} df_0 &=& f_1
f^{l+1}_{-1} 
df_0\\ &+&  (q^{-l}-1)f^l_{-1}df_0+l(1-q^{-1})f^l_{-1}df_0 ,\\
(f_0-1)f^l_{-1}df_1&=& q(f_0-1) f^{l-1}_{-1}df_0 \\
&-& q^{l+1} f_1(f_0-1) f^{l-1}_{-1} 
df_{-1}-q^2(1-q^{l-1})(f_0-1)f^{l-2}_{-1}df_{-1}, ~l>0.
\end{eqnarray*}
It follows that $\gamma \in ker \pi_{2_{\Gamma}}$ can be written in the form 
\begin{eqnarray*} 
\gamma&=&\sum_{k \geq 0} \tilde{a}_k f^k_1 (f_0-1) df_1
\\
&+& \sum_{k,l \geq 0} \tilde{c}_{k,l} (f_0-1) f^k_1 f^l_{-1} df_{-1}
+\sum_{k,l 
\geq 0} \tilde{g}_{k,l} f^k_1 f^l_{-1} df_0.
\end{eqnarray*}
Assuming $\gamma \in ker \pi_{1_{\Gamma}} \bigcap ker \pi_{2_{\Gamma}}$ one
obtains
\begin{eqnarray*} 
\pi_{1_{\Gamma}}(\gamma)&=& \sum_{k \geq 0} 
\tilde{a}_k x^k(xx^*-1) dx \\
&+& \sum_{k,l \geq 0} \tilde{c}_{k,l} (xx^*-1) x^k{x^*}^l dx^* + \sum_{k,l
\geq 
0} \tilde{g}_{k,l} x^k{x^*}^l d(xx^*)=0. 
\end{eqnarray*}
The left coefficient of $dx$ is
\[ \sum_{k \geq 0} (\tilde{a}_k x^{k+1} x^*- \tilde{a}_k x^k) + q^{-1}
\sum_{k,l 
\geq 0} q^{-1} \tilde{g}_{k,l} x^k {x^*}^{l+1} =0, \]
which gives
$ \tilde{a}_k=0,~\forall k$ and $\tilde{g}_{k,l}=0,~\forall k,l$. This
leads to 
\[ \pi_{1_\Gamma}(\gamma)=\sum_{k,l \geq 0} \tilde{c}_{k,l} (xx^*-1) 
x^k{x^*}^l=0. \]
Since $xx^*-1$ is not a zero divisor, it follows that $\tilde{c}_{k,l}=0~
\forall 
k,l$, i.e. $\gamma=0$.

Now, let us prove the assertion for $\Gamma^2(P(S^2_{qq\phi}))$. 
\\Applying $d$ to (\ref{refd1})-(\ref{refd4}) one obtains in
$\Gamma(P(S^2_{qq\phi}))$, besides
(\ref{refd21})-(\ref{refd24}), the following relations
\begin{eqnarray*}
 &&df_1 df_1=0,~~df_{-1} df_{-1}=0,~~ df_{-1} df_1 =-qdf_1
df_{-1},
\\&&df_0 df_1 =-df_1 df_0,~~df_{-1}df_0=-df_0 df_{-1},
\end{eqnarray*}
and one can see that $\Gamma^2(P(S^2_{qq\phi}))$ is generated by
the elements $df_1 df_{-1}$ and $df_0 df_0$ as a
$P(S^2_{qq\phi})$-bimodule.
Let us consider a general element $\gamma \in \Gamma^2(P(S^2_{qq\phi}))$,
\[ \gamma= a df_1 df_{-1} + b df_0 df_0. \]
Applying $\pi_{2_{\Gamma}}$ to $\gamma \in ker \pi_{2_{\Gamma}}$ gives
\[ \pi_2(a) dx dx^*=0, \] and it follows $a \in ker \pi_2$(Lemma \ref{modbas},
(ii)). $ker \pi_2$
is generated by the element $1-f_0$ (see (\ref{acoeff})), and because of
relation
(\ref{refd23}) $ker \pi_{2_{\Gamma}}$ is generated by $df_0df_0$. Now assume
\[\gamma = \tilde{b}df_0 df_0 \in ker \pi_{1_{\Gamma}} \bigcap ker
\pi_{2_{\Gamma}}. \] It follows that \begin{eqnarray*}
\pi_{1_\Gamma}(\gamma)=\pi_1(\tilde{b}) d(xx^*) d(xx^*) &
=&\pi_1(\tilde{b}) (x(dx^*)(dx)x^* + (dx)
x^*
x^*(dx)) \\ &=& \pi_1(\tilde{b})(x^* x dx dx^* + xx^*dx^*dx) \\
&=& \pi_1(\tilde{b})(x^*x dx dx^* -q xx^*dx dx^*)\\&=&
\pi_1(\tilde{b})(1-q)dx
dx^*=0, \end{eqnarray*} and this leads to $\tilde{b} \in ker \pi_1$.

As noted above $ker \pi_1$ is generated by the element $f_1f_{-1} -
f_0$.
It is immediate from (\ref{re1})-(\ref{re4}) that
$(f_1f_{-1}-f_0)f_i=q^if_i(f_1f_{-1}-f_0),~~i=\pm 1,0$. This and
(\ref{refd24})
gives $\gamma =0$, i.e. $(ker \pi_{1_{\Gamma}})^2 \bigcap
(ker\pi_{2_{\Gamma}})^2=0$.

Finally, $\Gamma^n(P(S^2_{qq\phi}))=0,~~\forall n>2$, i.e.
$J^n(P(S^2_{qq\phi}))=\Omega^n(P(S^2_{qq \phi})),~~\forall n>2$,
since $\Gamma^n(P(D_q))=0,~~\forall n>2$.
This can also be obtained directly from the generators of
$J(P(S^2_{qq\phi}))$
((\ref{refd1})- (\ref{refd4}) and (\ref{refd21})- (\ref{refd24})), and we
need
no additional relations for higher degrees. \hfill$\square$
\\\\$\Gamma(P(S^2_{qq\phi}))$ is a 
$*$-differential calculus:

Obviously, the homomorphisms $\pi_1$ and $\pi_2$ are $*$-homomorphisms. We
know
that the universal differential calculus $\Omega(P(S^2_{qq\phi}))$ is a
$*$-differential calculus, and one easily verifies that the homomorphisms
$\pi_{1_{\Omega \rightarrow \Gamma}}$ and $\pi_{2_{\Omega \rightarrow
\Gamma}}$
are $*$-homomorphisms. It follows that the kernels of these homomorphisms
are
$*$-ideals and also the intersection of these kernels is a $*$-ideal. This
ideal
is just the differential ideal defining the differential calculus
$\Gamma(P(S^2_{qq\phi}))$, i.e. there exists a $*$-structure on
$\Gamma(P(S^2_{qq\phi}))$, such that the quotient map
$\pi_{\Omega,\Gamma}: \Omega(P(S^2_{qq\phi})) \longrightarrow
\Gamma(P(S^2_{qq\phi}))$ satisfies \[ \pi_{\Omega,\Gamma} \circ *=* \circ
\pi_{\Omega,\Gamma}. \]
\\Finally, the differential calculus
$\Gamma(P(S^2_{qq\phi}))$ is also $U_{q^{1/2}}(sl_2)$-covariant:

First we recall the meaning of the covariance of the differential calculus
over
the quantum disc with respect to the action of $U_{q^{1/2}}(sl_2)$ and refer for more details to \cite{vaks}.

The Hopf algebra $U_{q^{1/2}}(sl_2)$ is the algebra generated by
$K^{\pm 1}$, $E$, $F$ with the relations \begin{eqnarray*}
&& KK^{-1}=K^{-1}K=1,~~K^{\pm 1}E=q^{\pm 1} E K^{\pm 1}, ~~K^{\pm 1}
F = q^{\mp 1} F K^{\pm 1}  \\
&&EF-FE=(K-K^{-1})/(q^{1/2}- q^{-1/2}). \end{eqnarray*} The comultiplication
$\Delta$, the counit $\varepsilon$ and the antipode $S$ are defined as
follows \begin{eqnarray*}
&& \Delta(E)=E \otimes 1 + K \otimes E, ~\Delta(F)= F\otimes K^{-1} + 1
\otimes
F,\\ && \Delta(K^{\pm 1}) = K^{\pm 1} \otimes K^{\pm 1}, \\
&& \varepsilon(E)= \varepsilon(F)=0,~~ \varepsilon(K^{\pm 1})=1, \\
&& S(E)=-K^{-1}E,~~S(F)=-FK,~~ S(K^{\pm 1})= K^{\mp 1}.
\end{eqnarray*}
There exists an action $\cdot: U_{q^{1/2}}(sl_2) \times P(D_q)
\longrightarrow
P(D_q)$, which means
\begin{eqnarray*} h \cdot 1&=&\varepsilon(h)1,~~\forall h \in
U_{q^{1/2}}(sl_2) \\
1 \cdot a &=&a,~~\forall a \in P(D_q) \\
h \cdot (ab) &=&\sum (h_1 \cdot a) (h_2 \cdot b),~~\forall h \in
U_{q^{1/2}}(sl_2),~~a,b \in P(D_q)\\
h \cdot(g \cdot a) &=& (hg) \cdot a,~~\forall h,g \in U_{q^{1/2}}(sl_2),~~a,
\in P(D_q), \end{eqnarray*} defined by
\begin{eqnarray*} && K^{\pm 1} \cdot x :=q^{ \pm 1} x,~~ F \cdot x:=
q^{1/4}1 ,~~ E \cdot x:=-q^{1/4} x^2 \\
&& K^{\pm} \cdot x^* := q^{\mp 1} x^*,~~F \cdot x^*:=- q^{5/4} {x^*}^2,~~E
\cdot x^* := q^{-3/4}1. \end{eqnarray*}

It was shown in \cite{vaks} that this action can be extended to the
differential
algebra $\Gamma(P(D_q))$, which means that there exists $\cdot:
U_{q^{1/2}}(sl_2) \times \Gamma(P(D_q)) \longrightarrow \Gamma(P(D_q))$
defined by \[ h \cdot da = d( h \cdot a), \]
and
\[ h \cdot \gamma \sigma = \sum (h_1 \cdot \gamma) ( h_2 \cdot \sigma),~~ h
\in U_{q^{1/2}}(sl_2),~~\gamma, \sigma \in \Gamma(P(D_q)). \]
It turns out that this action is compatible with the gluing procedure, i.e.
 there is an action $\cdot: U_{q^{1/2}}(sl_2) \times
P(S^2_{qq\phi}) \longrightarrow  P(S^2_{qq\phi})$ given by
\[ h \cdot (a,b):= (h \cdot a, h \cdot b),~~h \in U_{q^{1/2}}(sl_2)
,~~(a,b) \in P(S^2_{qq \phi}). \]
On the generators $f_1$, $f_{-1}$ and $f_0$ we have
\begin{eqnarray*} && K^{\pm 1} \cdot f_1= q^{\pm 1} f_1,~~F \cdot f_1 =
q^{1/4} 1,~~ E \cdot f_1 =-q^{1/4} f^2_1, \\
&&K^{\pm 1} \cdot f_{-1} = q^{\mp 1} f_{-1},~~F \cdot f_{-1} = -q^{5/4}
f^2_{-1},~~E \cdot f_{-1}= q^{-3/4} 1, \\
&& K^{\pm 1} \cdot f_0 =f_0,~~F \cdot f_0=q^{5/4}(f_{-1}-f_0f_{-1}),~~E
\cdot
f_0 = q^{1/4}(f_1- f_1f_0). \end{eqnarray*}

The extension to $\Gamma(P(S^2_{qq\phi}))$ is obvious. The projections
$\pi_1$
and $\pi_2$ intertwine the actions of $U_{q^{1/2}}(sl_2)$ on
$P(S^2_{qq\phi})$ and on the two copies of $P(D_q)$. This property extends
to
the universal differential calculus $\Omega(P(S^2_{qq\phi}))$, and we even
have \[
\pi_{i_{\Omega \rightarrow \Gamma}} (h \cdot \gamma)= h \cdot \pi_{i_{\Omega
\rightarrow \Gamma}}(\gamma), \] which means that
the kernels of the homomorphisms $\pi_{i_{\Omega \rightarrow \Gamma}}$ are
invariant under the action $\cdot$. Thus, also $\ker \pi_{1_{\Omega
\rightarrow \Gamma}} \bigcap \ker \pi_{2_{\Omega \rightarrow \Gamma}}$ is
invariant under the action $\cdot$. This intersection is
just our differential ideal $J(P(S^2_{qq\phi}))$, and it follows that one
can
extend the action $\cdot$ on $P(S^2_{qq\phi})$ to the differential
calculus
$\Gamma(P(S^2_{qq\phi}))$, i.e. $\Gamma(P(S^2_{qq\phi}))$ is covariant.

It is also easy to show that $P(S^2_{qq\phi})$ and $\Gamma(P(S^2_{qq\phi}))$ 
are covariant with respect to $U_{q^{1/2}}(su(1,1))$ as $*$-algebras. 
(cf.\cite{vaks})

\section{Appendix}
\indent

The purpose of this appendix is to prove Lemma \ref{modbas}.\\
(i) We have to show that from $ adx+bdx^*=0,~a,b \in
P(D_q)$
follows $a=0$ and $b=0$. Recall that
\[\Omega^1(P(D_q))=\{\sum_k a_k \otimes b_k| \sum_k a_k b_k=0\} \subset
P(D_q)
\otimes P(D_q) \]
$(da =1 \otimes a - a \otimes 1)$.
By Lemma \ref{basis} we have a left $P(D_q)$-module basis in $P(D_q) \otimes
P(D_q)$
formed by the elements $1 \otimes x^k{x^*}^l$,
therefore the elements $\{d(x^k {x^*}^l),~k,l >0\}$ form a left
$P(D_q)$-module
basis in
$\Omega^1(P(D_q))$: From $\sum_{kl} a_{kl}d(x^k {x^*}^l)=
\sum_{kl} a_{kl} \otimes x^k {x^*}^l - \sum a_{kl} x^k
{x^*}^l\otimes 1=0,~a_{kl}
\in P(D_q),$ 
follows $a_{kl}=0,~\forall k,l$.

Now we define the following left module projection $P_1 : 
\Omega^1(P(D_q))
\longrightarrow \{adx+bdx^*|a,b\in P(D_q)\}\subset
\Omega^1(P(D_q))$: \begin{eqnarray*}
P_1(dx) &:=& dx,~~P_1(dx^*):=dx^* \\
P_1(d(x^k)) &:=& \sum^{k-1}_{i=0} q^i x^{k-1} dx,~~ k>0 \\
P_1(d({x^*}^l))&:=& \sum^{l-1}_{i=0} q^{-i} {x^*}^{l-1} dx^*,~~l>0 \\
P_1(d(x^k{x^*}^l))&:=& \sum^{k-1}_{i=0} q^{i-l} x^{k-1} {x^*}^ldx +
\sum^{l-1}_{i=0} q^{-i} x^k {x^*}^{l-1} dx^*,~~ k>0;~l>0  \end{eqnarray*}

All generators of $J(D_q)$ lie in the kernel of $P_1$
(for example: $P_1(xdx-q^{-1}(dx)x= P_1(q^{-1}((1+q)xdx-dx^2)=0$), thus the
left
module generated by these elements lies in the kernel of $P_1$. If we can
show
that also the right module generated by these elements lies in the kernel of
$P_1$,
the whole first degree of $J(P(D_q))$ lies in $ker P_1$. We show this for
the generator $xdx -q^{-1}(dx)x=(1+q^{-1})xdx-q^{-1}dx^2$ and leave the
remaining cases to the reader.

The Leibniz rule gives
\[((1+q^{-1})xdx-q^{-1}dx^2) x^k {x^*}^l\\
=(1+q^{-1})xd(x^{k+1}{x^*}^l) -x^2d(x^k
{x^*}^l)- q^{-1}d(x^{k+2}{x^*}^l). \] Applying $P_1$ to this formula one
gets
\begin{eqnarray*} && P_1(((1+q^{-1})xdx -q^{-1}dx^2)x^k{x^*}^l)\\ &=&
(1+q^{-1})
\sum^k_{i=0} q^{i-l} x^{k+1}{x^*}^ldx - (1+q^{-1})\sum^{l-1}_{i=0} q^{-i}
x^{k+2}{x^*}^{l-1}dx^* \\
&-& \sum^{k-1}_{i=0} q^{i-l} x^{k+1}{x^*}^ldx -\sum^{l-1}_{i=0}q^{-i}
x^{k+2}
{x^*}^{l-1} dx^* \\
&-& q^{-1} \sum^{k+1}_{i=0} q^{i-l} x^{k+1} {x^*}^ldx - q^{-1}
\sum^{l-1}_{i=0}
q^{-i}x^{k+2} {x^*}^{l-1} dx =0,~~ \forall k,l. \end{eqnarray*}

The calculation for the remaining generators of $J(P(D_q))$ is analogous,
thus
\begin{equation} \label{prop1} J^1(P(D_q)) \subseteq ker P_1.
\end{equation}

Let $\pi_{\Omega,\Gamma} : \Omega(P(D_q)) \longrightarrow \Gamma(P(D_q))$
be the quotient map. Because of
formula (\ref{prop1}) there exists a left $P(D_q)$-module homomorphism
$\Lambda_1:\Gamma^1(P(D_q)) \longrightarrow \Omega^1(P(D_q))$ defined by
\[ \Lambda_1 \circ \pi_{\Omega,\Gamma}:=P_1. \] Applying $\Lambda_1$ to $adx
+ bdx^*=0$ it follows that
$adx+bdx^*=0$ in $\Omega^1(P(D_q))$, which gives $a=0$ and $b=0$. The right
basis property follows now easily from the left one and the relations 
defining $\Gamma(P(D_q))$.\\
(ii) Applying the
differential $d$ to the generators of $J(P(D_q))$ one obtains the
following relations in $\Gamma^2(P(D_q))$:
\[ dxdx=0;~~dx^*dx^*=0;~~dx^* dx=-qdx dx^*. \] 

The procedure of the proof of (i) carries over to the
degree two case:

First one proves that the elements $d(x^m{x^*}^n)d(x^k{x^*}^l)$ form a
left
$P(D_q)$-module basis of $\Omega^2(P(D_q))$.
Then one defines a left $P(D_q)$-module projection\\
$ P_2: \Omega^2(P(D_q))
\longrightarrow \{adx dx^*,~a \in P(D_q)\} \subset \Omega^2(P(D_q))$:
\begin{eqnarray*}
&& P_2(dx dx^*):=dx dx^*,~~P_2(dx^* dx):=-qdx dx^* \\
&& P_2(dx^m dx^k):=0,~~P_2(d{x^*}^n d{ x^*}^l):=0 \\
&& P_2(dx^m d{x^*}^l):= q^{1-l} \sum^{m-1}_{i=0} \sum^{l-1}_{s=0} q^{i-s}
x^{m-1} {x^*}^{l-1} dx dx^*,~~m,l \geq 1 \\
&& P_2(d{x^*}^n dx^k):=-q^k \sum^{n-1}_{i=0} \sum^{k-1}_{s=0} q^{s-i}
{x^*}^{n-1} x^{k-1} dx dx^*,~~n,k \geq 1 \\
&& P_2(d(x^m {x^*}^n) d(x^k {x^*}^l):=q^{k-l}\sum^{m-1}_{i=0}
\sum^{l-1}_{s=0}
q^{i-n} q^{1-s} x^{m-1} {x^*}^n x^k {x^*}^{l-1}dx dx^* \\
&-& q^{k-l}\sum^{n-1}_{i=0}
\sum^{k-1}_{s=0} q^{-i} q^{s-l} x^m{x^*}^{n-1} x^{k-1} {x^*}^ldx
dx^*,~~m,n,k,l
\geq 1. \end{eqnarray*}

With these definitions one proves that $J^2(P(D_q)) \subseteq ker P_2$.

Thus, there exists a left $P(D_q)$-module homomorphism
$\Lambda_2: \Gamma^2(P(D_q)) \longrightarrow
\Omega^2(P(D_q))$ defined by \[\Lambda_2 \circ
\pi_{\Omega,\Gamma}:=P_2;~~(\Lambda_2(adxdx^*)=adxdx^*) \]
Applying $\Lambda_2$ to $adx dx^*=0,~~a \in P(D_q)$ in $\Gamma(P(D_q))$ it
follows  that $adxdx^*=0$ in $\Omega^2(P(D_q))$, which gives $a=0$. Again, the
right basis property is now easily derived.\\
(iii) Immediate.
\hfill$\square$
\\\\ {\bf Acknowledgments:} We thank P. M. Alberti and P. M. Hajac for helpful 
discussions and P. Podle\'s and W. Pusz for valuable remarks.

\end{document}